\newif\ifprivate\privatefalse \newif\iflongversion\longversionfalse \newif\ifoldversion\oldversionfalse \def\loadTIKZ{\usepackage{tikz}\usetikzlibrary{matrix,arrows,calc,cd,decorations.pathmorphing}}\ifdefined\headpresent\else \providecommand\DocumentClassOptions{article}\ifdefined\presentationON 
\documentclass[\DocumentClassOptions,20pt]{amsart}\usepackage{extsizes} \else \documentclass[\DocumentClassOptions]{amsart}\fi \usepackage{ifpdf}\ifdefined\BibLaTeX \usepackage[\BibLaTeX]{biblatex}\DeclareFieldFormat{labelalpha}{\thefield{entrykey}}\makeatletter\protected\def\abx@missing#1{\mbox{\reset@font\color{red}#1??}}\makeatother \DeclareFieldFormat{title}{\textit{#1}}\DeclareFieldFormat{year}{\textbf{#1}}\renewbibmacro*{date}{\printtext[year]{\printdate}}\BibLaTeXCommands \else\fi \newcommand\DIRroot{../../../../../../../../../../../}\overfullrule=5pt \newcommand\hfuzzReset{\hfuzz=3pt}\hfuzzReset \newcommand\toleranceReset{\tolerance=1400}\toleranceReset \newcommand\emergencystretchReset{\emergencystretch=1ex}\emergencystretchReset \hbadness=10000 \usepackage{xifthen}\usepackage{forarray}\usepackage{xstring}\usepackage{stringstrings}\def\StackCreate#1#2#3{\expandafter\def\csname#1\endcsname{#3}\expandafter\def\csname#1Push\endcsname##1{\expandafter\edef\csname#1\endcsname{##1#2\csname#1\endcsname}}\expandafter\def\csname TopAux#1\endcsname ##1#2##2#3{##1}\expandafter\def\csname#1Top\endcsname{\expandafter\expandafter\expandafter\expandafter\expandafter\expandafter\csname TopAux#1\endcsname\csname#1\endcsname}\expandafter\def\csname PopAux#1\endcsname ##1#2##2#3##3#2{\expandafter\def\csname##3\endcsname{##2#3}}\expandafter\def\csname#1Pop\endcsname{\expandafter\expandafter\expandafter\expandafter\expandafter\expandafter\csname PopAux#1\endcsname\csname#1\endcsname#1#2}}\def\GetAfterColonAux#1:#2;{#2}\def\GetAfterColon#1{\IfBeginWith{#1}{:}{\GetAfterColonAux#1;}{#1}}\usepackage{aliascnt}\newlength{\tempwidth}\newcommand{\fillX}[2][]{\settowidth{\tempwidth}{#2}\def\temp{#1}\ifx\temp\empty\else\addtolength{\tempwidth}{#1}\fi\leavevmode\cleaders\hbox to \tempwidth{\hss #2\hss }\hfill\kern0pt }\usepackage[shortlabels,inline]{enumitem}\setenumerate[1]{leftmargin=5.5ex}\setitemize[1]{leftmargin=5.5ex}\SetEnumitemKey{noindent}{leftmargin=0ex, itemindent=5ex, align=right, itemsep=1ex }\newcounter{ManualLabel}\makeatletter \newcommand\itemPatch[1][]{\item[\theenumi#1]\refstepcounter{ManualLabel}\def\@currentlabel{\theenumi#1}}\newcommand\itemDescribe[1][]{\item[#1]\refstepcounter{ManualLabel}\def\@currentlabel{#1}}\makeatother \ifdefined\selectPages \usepackage[\selectPages]{pagesel} \fi \ifdefined\presentationON \usepackage[hcentering=true,textheight=240mm,textwidth=210mm]{geometry} \else\fi \usepackage{everypage-1x}\ifdefined\AddEverypageHook \newcommand\AddPrivateToMargin[1]{\AddEverypageHook{\tikz[overlay,remember picture]{\node at ($(current page.west)+(1.5,0)$) [rotate=90] {\textcolor{orange}{\vbox{\hrule width \the\textwidth height 0.5pt} \textcolor{defaultcolor}{#1}\ \vbox{\hrule width 40em height 0.5pt}}}; }}}\newcommand\AddLongversionToMargin[1]{\AddEverypageHook{\tikz[overlay,remember picture]{\node at ($(current page.west)+(2,0)$) [rotate=90] {\textcolor{\LongColor}{\vbox{\hrule width \the\textwidth height 0.5pt} \textcolor{defaultcolor}{#1}\ \vbox{\hrule width 40em height 0.5pt}}}; }}}\newcommand\AddOldversionToMargin[1]{\AddEverypageHook{\tikz[overlay,remember picture]{\node at ($(current page.west)+(2.5,0)$) [rotate=90] {\textcolor{\OldColor}{\vbox{\hrule width \the\textwidth height 0.5pt} \textcolor{defaultcolor}{#1}\ \vbox{\hrule width 40em height 0.5pt}}}; }}}\newcommand\AddLineToMargin[3]{\AddEverypageHook{\tikz[overlay,remember picture]{\node at ($(current page.west)+(#2,0)$) [rotate=90] {\textcolor{#1}{\vbox{\hrule width \the\textwidth height 0.5pt} \textcolor{defaultcolor}{#3}\ \vbox{\hrule width 40em height 0.5pt}}}; }}}\else\fi \providecommand\AddPrivateToMargin[1]{\AddToHook {shipout/background}{\tikz[overlay,remember picture]{\node at ($(current page.west)+(1.5,0)$) [rotate=90] {\textcolor{orange}{\vbox{\hrule width \the\textwidth height 0.5pt} \textcolor{black}{#1}\ \vbox{\hrule width 40em height 0.5pt}}}; }}}\providecommand\AddLongversionToMargin[1]{\AddToHook {shipout/background}{\tikz[overlay,remember picture]{\node at ($(current page.west)+(2,0)$) [rotate=90] {\textcolor{\LongColor}{\vbox{\hrule width \the\textwidth height 0.5pt} \textcolor{black}{#1}\ \vbox{\hrule width 40em height 0.5pt}}}; }}}\providecommand\AddOldversionToMargin[1]{\AddToHook {shipout/background}{\tikz[overlay,remember picture]{\node at ($(current page.west)+(2.5,0)$) [rotate=90] {\textcolor{\OldColor}{\vbox{\hrule width \the\textwidth height 0.5pt} \textcolor{black}{#1}\ \vbox{\hrule width 40em height 0.5pt}}}; }}}\providecommand\AddLineToMargin[3]{\AddToHook {shipout/background}{\tikz[overlay,remember picture]{\node at ($(current page.west)+(#2,0)$) [rotate=90] {\textcolor{#1}{\vbox{\hrule width \the\textwidth height 0.5pt} \textcolor{black}{#3}\ \vbox{\hrule width 40em height 0.5pt}}}; }}}\let\oldtocsection=\tocsection \let\oldtocsubsection=\tocsubsection \let\oldtocsubsubsection=\tocsubsubsection \renewcommand{\tocsection}[2]{\hspace{0em}\vspace*{0.1em}\oldtocsection{#1}{#2}}\renewcommand{\tocsubsection}[2]{\hspace{4ex}\oldtocsubsection{#1}{#2}}\renewcommand{\tocsubsubsection}[2]{\hspace{6ex}\oldtocsubsubsection{#1}{#2}}\usepackage{ulem}\usepackage{fancybox}\ifpdf \usepackage[pdftex]{lscape}\else \usepackage{lscape}\fi \makeatletter \newcommand{\verbatimfont}[1]{\def\verbatim@font{#1}}\makeatother \IfFileExists{mathabx.sty}{}{}\usepackage{amsfonts}\usepackage{amssymb}\usepackage{stmaryrd}\usepackage{amsmath}\usepackage{amsthm}\usepackage{dsfont}\IfFileExists{mbboard.sty}{\usepackage{mbboard}}{}\usepackage{mathrsfs}\usepackage{twcal}\usepackage{accents}\usepackage{bm}\usepackage[T1]{fontenc}\usepackage[latin1]{inputenc}\catcode`\=13 \def{+}\newcommand\assigncharacter[1]{\expandafter\newcommand\csname #1\endcsname{\mathds{#1}}}\FunctionForEach{,}{\assigncharacter}{A,B,C,D,E,F,G,I,J,K,M,N,Q,R,T,U,V,W,X,Y,Z}\renewcommand\assigncharacter[1]{\expandafter\newcommand\csname C#1\endcsname{\mathcal{#1}}}\FunctionForEach{,}{\assigncharacter}{A,B,C,D,E,F,G,H,I,J,K,L,M,N,O,P,Q,R,S,T,U,V,W,X,Y,Z}\renewcommand\assigncharacter[1]{\expandafter\newcommand\csname D#1\endcsname{\mathfrak{#1}}}\FunctionForEach{,}{\assigncharacter}{a,b,c,d,e,f,g,h,i,j,k,l,m,n,o,p,q,r,s,t,u,v,w,x,y,z,A,B,C,D,E,F,G,I,K,L,M,N,O,P,Q,R,S,T,U,V,W,X,Y,Z} \renewcommand\assigncharacter[1]{\expandafter\newcommand\csname S#1\endcsname{\mathscr{#1}}}\FunctionForEach{,}{\assigncharacter}{A,B,C,D,E,F,G,H,I,J,K,L,M,N,O,P,Q,R,T,U,V,W,X,Y,Z}\def\NewFont#1#2#3#4#5{\expandafter\font\csname #1display\endcsname =#1 at #2 \expandafter\font\csname #1normal\endcsname =#1 at #3 \expandafter\font\csname #1script\endcsname =#1 at #4 \expandafter\font\csname #1scriptscript\endcsname =#1 at #5 }\def\NewFontLetter#1#2{{\mathchoice {{\expandafter\hbox{\csname #1display\endcsname\char"#2}}}{{\expandafter\hbox{\csname #1normal\endcsname\char"#2}}}{{\expandafter\hbox{\csname #1script\endcsname\char"#2}}}{{\expandafter\hbox{\csname #1scriptscript\endcsname\char"#2}}}}}\NewFont{pxsyc}{9.00pt}{8.00pt}{7.00pt}{6.00pt}\NewFont{pxsya}{9.00pt}{8.00pt}{7.00pt}{6.00pt}\renewcommand{\rightsquigarrow}{\NewFontLetter{pxsya}{20}}\NewFont{p1xr}{10.00pt}{9.00pt}{8.00pt}{7.00pt}\NewFont{MnSymbolA5}{10.00pt}{9.00pt}{8.00pt}{7.00pt}\NewFont{MnSymbolC10}{10.00pt}{9.00pt}{8.00pt}{7.00pt}\NewFont{MnSymbolD10}{12.00pt}{11.00pt}{10.00pt}{9.00pt}\NewFont{MnSymbolF10}{12.00pt}{11.00pt}{10.00pt}{9.00pt}\renewcommand{\bigcup}{\mathop{\NewFontLetter{MnSymbolF10}{1F}}}\def\IndependenceX#1#2{#1\setbox0=\hbox{$#1x$}\kern\wd0\hbox to 0pt{\hss$#1\mid$\hss}\lower.9\ht0\hbox to 0pt{\hss$#1\smile$\hss}\kern\wd0}\def\nIndependenceX#1#2{#1\setbox0=\hbox{$#1x$}\kern\wd0 \hbox to 0pt{\mathchardef\nn=12854\hss$#1\nn$\kern1.4\wd0\hss}\hbox to 0pt{\hss$#1\mid$\hss}\lower.9\ht0 \hbox to 0pt{\hss$#1\smile$\hss}\kern\wd0}\NewFont{manfnt}{12.00pt}{11.00pt}{10.00pt}{9.00pt}\NewFont{favmr7y}{12.00pt}{11.00pt}{10.00pt}{9.00pt}\ifpdf \usepackage[pdftex,usenames,x11names]{xcolor}\else \usepackage[dvips,usenames,x11names]{xcolor}\fi \StackCreate{ColoR}{;}{?}\AtBeginDocument{\colorlet{defaultcolor}{.}\ColoRPush{defaultcolor}}\definecolor{Green}{rgb}{0.00,0.50,0.00}\definecolor{DarkGreen}{rgb}{0.00,0.40,0.00}\definecolor{gray}{rgb}{0.40,0.40,0.40} \renewcommand\textcolor[2]{\ColoRPush{#1}\color{\ColoRTop}#2\ColoRPop\color{\ColoRTop}}\usepackage[pdftex]{graphicx}\usepackage[all]{xy}\ifdefined\loadTIKZ \loadTIKZ \def\TIKZlabel#1{}\else \usepackage{tikz}\usetikzlibrary{matrix,arrows,calc,cd,decorations.pathmorphing,intersections}\fi \newcommand {\notion}[2][]{\def\temp{#1}\ifmmode #2 \ifx \temp\empty \index{$#2$}\else \index{$#1$}\fi \else {\bf #2}\ifx \temp\empty \index{#2}\else \index{#1}\fi \fi }\newcommand\NOPAGENUMBER[1]{}\usepackage{xr-hyper}\newcommand{\refX}[2]{\IfBeginWith{#1}{:}{\ref{\GetAfterColonAux#1;-#2}}{\cite[\ref{#1-#2}]{#1}}}\providecommand\hyperrefBOOKMARKS{true}\ifdefined\BibLaTeX \usepackage[destlabel=true,pdftex,linktocpage,breaklinks,bookmarks=\hyperrefBOOKMARKS]{hyperref}\else \usepackage[destlabel=true,pdftex,linktocpage,pagebackref,breaklinks,bookmarks=\hyperrefBOOKMARKS]{hyperref}\fi \hypersetup{bookmarksdepth=3, colorlinks=true,allcolors=Green, linkcolor=DarkGreen, citecolor=violet, urlcolor=blue, runcolor=red, filecolor=magenta }\providecommand\url[1]{}\providecommand\nolinkurl[1]{}\providecommand\href[3][]{}\providecommand\hyperlink[2]{}\providecommand\hypertarget[2]{}\providecommand\hyperdef[3]{}\providecommand\hyperref[2][]{} \providecommand\hypersetup[1]{}\providecommand\pdfbookmark[3][]{}\providecommand\currentpdfbookmark[2]{}\providecommand\belowpdfbookmark[2]{}\providecommand\texorpdfstring[2]{} \def\UndefinedRef#1{\LARGE\bfseries\color{red} ??#1??}\makeatletter \def\@setref#1#2#3{\ifx#1\relax \protect\G@refundefinedtrue \nfss@text{\reset@font\UndefinedRef{#3}}\@latex@warning{Reference `#3' on page \thepage \space undefined }\else \expandafter\Hy@setref@link#1\@empty\@empty\@nil{#2}\fi }\makeatother \fi \newcommand\pr{\begin{proof}}\def\ende{\end{proof}}\newtheoremstyle{LayoutVoid}{1ex}{0ex}{\normalfont}{}{\bf}{.}{1ex}{}\newcommand\stressstatement[1]{#1}\theoremstyle{plain}\swapnumbers \newcommand\maketheorem[1]{\newtheorem{#1}[theorem]{\stressstatement{#1}} \newtheorem{#1Definition}[theorem]{\stressstatement{#1 and Definition}}  }\FunctionForEach{,}{\maketheorem}{Conclusion,Conjecture,Corollary,Fact,Facts,Lemma,Observation,Observations,Proposition,Reminder,Scholium,Summary,Theorem}\theoremstyle{definition}\theoremstyle{remark}\FunctionForEach{,}{\maketheorem}{Convention,Counterexample,Counterexamples,Discussion,Example,Examples, Exercise,Exercises,Explanation,Notation,Project,Projects,Question,Questions,Remark,Remarks,Strategy,Warning}\theoremstyle{LayoutVoid}\numberwithin{equation}{section}\newcommand{\labelon}[1]{\marginpar{#1}}\newcommand{\labelx}[1]{{\def\temp{#1}\ifx\temp\empty\else \label{#1}\labelon{#1}\fi}}\def\GetAfterColon#1:#2;;{#2}\def\GetAfterPlus#1+#2;;{#2}\newenvironment{FACT}[2]{\edef\currentlabel{#2}\IfBeginWith{#1}{:}{\def\tempFactName{void}\def\tempFreeTitle{\GetAfterColon#1;;\ }}{\IfBeginWith{#1}{+}{\def\tempFactName{voidTheorem}\def\tempFreeTitle{\GetAfterPlus#1;;\ }}{\def\tempFactName{#1}\def\tempFreeTitle{}}}\def\tempfacT{\end{\tempFactName}}\begin{\tempFactName}\labelx{#2}\textup{\textbf{\tempFreeTitle}}\capitalize[q]{#1}\caselower[q]{#1}\global\edef\factname{\thestring}}{\tempfacT}\newcommand{\dd}{\mathrm{d}}\catcode95=12 \catcode95=8 \newcommand\bdl{{\ifmmode \mathrm{bdlat}\else {bounded distributive lattice}\fi}} \newcommand{\st}{{\ \vert \ }}\let\temp\phi \let\phi\varphi \let \varphi\temp \let\temp\theta \let\theta\vartheta \let \vartheta\temp \let\eps\varepsilon  \let\0\emptyset \newcommand{\into}{\hookrightarrow}\makeatletter \newcommand{\xRightarrow}[2][]{\ext@arrow 0359\Rightarrowfill@{#1}{#2}}\newcommand{\xLeftarrow}[2][]{\ext@arrow 0359\Leftarrowfill@{#1}{#2}}\newcommand{\xonto}[2][]{\ext@arrow 0359\rightarrowfill@ {#1}{#2}\mathrel{\mspace{-15mu}}\rightarrow}\newcommand{\xinto}[2][]{\lhook\joinrel\ext@arrow 0359\rightarrowfill@ {#1}{#2}}\makeatother \newcommand{\lra}{\longrightarrow}\newcommand{\Ra}{\Rightarrow}\newcommand{\La}{\Leftarrow}\renewcommand{\>}{\rangle}\newcommand{\spez}{\,\rightsquigarrow\,}\newcommand{\alg}{\mathrm{alg}}\newcommand{\mal}{{\cdot}} \newcommand{\emailTressl}{marcus.tressl@manchester.ac.uk}\newcommand{\homepageTressl}{\url{http://personalpages.manchester.ac.uk/staff/Marcus.Tressl/}}\newcommand{\monthname}[1]{\ifcase#1 \or January \or February \or March \or April \or May \or June \or July \or August \or September \or October \or November \or December \fi}\newcommand\PrivateColor{Brown4}\newcommand\LongColor{teal}\newcommand\OldColor{gray}\newcommand\PrivateStart{\ColoRPush{\PrivateColor}\color{\ColoRTop}}\newcommand\PrivateEnd{\ColoRPop\color{\ColoRTop}}\newcommand\LongStart{\ColoRPush{\LongColor}\color{\ColoRTop}[BEGIN LONG VERSION]}\newcommand\LongEnd{[END LONG VERSION]\ColoRPop\color{\ColoRTop}}\newcommand\COL{\ifmmode\colon\else :\ \fi}\newcommand{\claim}{\textit{Claim.}\ }\newcommand{\case}[1]{\textit{Case #1.}\ } \renewcommand{\mod}{{\operator{\,mod\,}}}\newcommand\operator[1]{\mathop{\operatorname{#1}}\nolimits}\newcommand\kon{{\mathcal{K}}}\newcommand\qcop{\mathop {\mathring{\kon}}\nolimits} \newcommand{\DCF}{{\rm DCF}}\newcommand{\CODF}{{\rm CODF}}\newcommand{\Quot}{\operator{Quot}}\newcommand{\card}{\operator{card}}\newcommand{\qf}{\operator{qf}}\newcommand{\ord}{\operator{ord}}\newcommand{\Spec}{\operator{Spec}}\newcommand{\Sped}{\operator{Sped}}\newcommand{\trdeg}{\operator{tr.deg}}\IfFileExists{C:/wb/System64/WinBatch.exe}{}{}\ifdefined\isinput\endinput\else\fi \let\temp\theta \let\theta\vartheta \let \vartheta\temp \makeatletter \def\LST#1{\def\tempa{#1}\futurelet\next\CreateCite@check}\def\CreateCite@check{\ifx\next\bgroup\expandafter\CreateCite@ii\else\expandafter\CreateCite@i\fi}\def\CreateCite@ii#1{\cite[{\ref{LeSTre2024-\tempa}\ref{LeSTre2024-#1}}]{LeSTre2024}}\def\CreateCite@i{\cite[{\ref{LeSTre2024-\tempa}}]{LeSTre2024}}\makeatother \ifprivate \else \definecolor{LinkColor}{rgb}{0.00,0.20,0.00}\hypersetup{bookmarksdepth=3, colorlinks=true,allcolors=black, linkcolor=LinkColor, citecolor=violet, urlcolor=blue, runcolor=black, filecolor=black }\fi \def\Ind{\setbox0=\hbox{$x$}\kern\wd0\hbox to 0pt{\hss$\mid$\hss} \lower.9\ht0\hbox to 0pt{\hss$\smile$\hss}\kern\wd0}\def\Notind{\setbox0=\hbox{$x$}\kern\wd0\hbox to 0pt{\mathchardef \nn=12854\hss$\nn$\kern1.4\wd0\hss}\hbox to 0pt{\hss$\mid$\hss}\lower.9\ht0 \hbox to 0pt{\hss$\smile$\hss}\kern\wd0}\def \UC {\operatorname{UC}}\newcommand \ec {e.c.}\renewcommand{\labelon}[1]{}\renewcommand\LongStart{\ColoRPush{\LongColor}\color{\ColoRTop}}\renewcommand\LongEnd{\ColoRPop\color{\ColoRTop}}\ifprivate \AddPrivateToMargin{private version} \fi \iflongversion \AddLongversionToMargin{long version} \fi \ifoldversion \AddOldversionToMargin{old version included} \fi 
\begin{document} \title{On ordinary differentially large fields} \author{Omar Le\'on S\'anchez} \address{Omar Le\'on S\'anchez, The University of Manchester\\ Department of Mathematics\\ Oxford Road \\ Manchester, M13 9PL, UK} \email{omar.sanchez@manchester.ac.uk} \author{Marcus Tressl} \address{Marcus Tressl, The University of Manchester\\ Department of Mathematics\\ Oxford Road \\ Manchester, M13 9PL, UK \newline Homepage: \homepageTressl} \email{\emailTressl} \date{\today}  \thanks{\emph{Acknowledgments.} The first author was partially supported by EPSRC grant EP/V03619X/1} \begin{abstract} We provide a characterisation of differentially large fields in arbitrary characteristic and a single derivation in the spirit of Blum axioms for differentially closed fields. In the case of characteristic zero, we use these axioms to characterise differential largeness in terms of being existentially closed in the differential algebraic Laurent series ring, and we prove that any large field of infinite transcendence degree can be expanded to a differentially large field even under certain prescribed constant fields. As an application, we show that the theory of proper dense pairs of models of a complete and modelcomplete theory of large fields, is a complete theory. As a further consequence of the expansion result we show that there is no real closed and differential field that has a prime model extension in closed ordered differential fields, unless it is itself a closed ordered differential field. \end{abstract} \maketitle \setcounter{tocdepth}{1} \tableofcontents \section{Introduction} \noindent The class of differentially large fields was introduced and studied by the authors in \cite{LeSTre2024}. Evidenced by the results in that paper, this class can be considered as the differential analogue of the class of large fields. We recall that a field $K$ is said to be large (aka ample) if every irreducible algebraic variety over $K$ with a smooth $K$-point has a Zariski dense set of $K$-points. Equivalently, $K$ is \ec\ (existentially closed) in the field of formal Laurent series $K((t))$. \par \medskip \par In \cite{LeSTre2024}, a differential field $(K,\Delta)$ of characteristic zero with commuting derivations $\Delta=\{\delta_1,\dots,\delta_m\}$ is defined to be differentially large if $K$ is large (as a field) and for every differential extension $(L,\Delta)$, if $K$ is \ec\ in $L$ (as a field), then $(K,\Delta)$ is \ec\ in $(L,\Delta)$ as a differential field (see \LST{ECbasic} for the algebraic meaning of \ec). Several foundational properties and applications are explored in \cite{LeSTre2024}. In particular, it is shown that $(K,\Delta)$ is differentially large if and only if $(K,\Delta)$ is \ec\ in $(K((t_1,\dots,t_m)),\Delta)$, where the derivations on $K((t_1,\dots,t_m))$ are the natural ones extending those on $K$ that commute with meaningful sums and satisfy $\delta_j(t_i)=\frac{dt_i}{dt_j}$. \par \medskip \par The first-order characterisation of differential largeness provided in \LST{DLisUC} makes reference to the somewhat elaborate axiom scheme $\UC$ from \cite[4.5]{Tressl2005}. In \ref{DLAxiomsOrdinary} below we give a significant simplification of this axiom scheme in the ordinary case, i.e. the case of a single derivation, so $\Delta=\{\delta\}$. The new scheme resembles the Blum axioms for differentially closed fields of characteristic 0 ($\DCF_0$) and at the same time allows an extension of the notion of differential largeness to arbitrary characteristic (cf. \ref{defdifflarge}). In subsequent sections we give applications of our new simple description of differential largeness in the ordinary case as follows. \par \smallskip Henceforth we restrict to a single derivation. An immediate consequence of the new axioms is the new characterisation \ref{CODFnewAxioms} of closed ordered differential fields (\CODF), in the sense of Singer \cite{Singer1978a}, which does not make reference to the order. A further corollary (\ref{DiffLargeDvariety}) provides geometric axioms for differentially large fields in arbitrary characteristic in terms of D-varieties, in the spirit of the Pierce-Pillay axioms for $\DCF_0$, see \cite{PiePil1998}. \par \smallskip In the rest of the paper we readopt the characteristic zero assumption. In Section~\ref{powerseries}, we prove that differential largeness can be characterised in terms of being existentially closed in the differential algebraic formal Laurent series, see \ref{ECinDiffAlgSeries}. Our proof uses an approximation-type statement that resembles that of Denef-Lipshitz in \cite{DenLip1984}. We then use this to produce a new way (or rather an improvement of the construction in \cite{LeSTre2024} for the ordinary case) to construct differentially large fields using iterated differential algebraic Laurent series, see \ref{ConstructionDiffAlgPowerSeries}. \par \smallskip In section \ref{expansionslarge}, we show that for any ordinary differential field $(K,d)$ and any given large field $L\supseteq K$ of transcendence degree over $K$ at least the size of $K$, there is an extension $\delta$ of $d$ to $L$ such that $(L,\delta)$ is differentially large, see Theorem \ref{DLexistence}. This has two consequences: Firstly, large fields of infinite transcendence degree (over $\Q$) are characterized in \ref{AllLargeSupportDL} as exactly those fields that possess a derivation $d$ for which $(L,d)$ is differentially large (significantly generalizing an earlier result by Christian Michaux saying that $\R$ carries a \CODF\ structure). Secondly, we show in \ref{CODFprimModels} that no real closed field equipped with any derivation has a prime model extension in \CODF, unless it is already a \CODF; this strengthens a result from \cite{Singer1978a} stating that the theory \CODF\ does not have a prime model. \par Theorem \ref{DLexistence} is significantly strengthened in section \ref{SectionExpansionsWithConstants} in the case when the constant field $C$ of $K$ is dense for the \'etale open topology of $L$ (see \ref{DefnEtOpenOnK} for its definition). Namely we show in Theorem \ref{ExtendToDL} that $L$ can be expanded to a differentially large field whose constant field is algebraic over $C$. This Theorem has an interesting consequence for dense pairs of large fields: In Corollary \ref{ExamplesCompleteTheoriesOfPairs} we show that for any complete and model complete theory of large fields of characteristic 0 in the language $\SL$ of rings, possibly extended by constants, the theory of proper dense pairs of models of $T$ is complete. This extends A. Robinson's theorem saying that the theory of proper pairs of dense real closed fields is complete. \par \medskip\noindent By a \textit{differential ring} in this paper we always mean a commutative unital ring furnished with a single derivation. \section{Blum-style axioms for ordinary differentially large fields}\label{BlumAxioms} \noindent In \cite{LeSTre2024}, differentially large fields in characteristic zero were introduced. The definition there makes sense also for ordinary differential fields of characteristic $p>0$.\looseness=-1 \begin{FACT}{Definition}{defdifflarge} A differential field $(K,d)$, of arbitrary characteristic, is said to be differentially large if it is large as a field and for every differential field extension $(L,\delta)/(K,d)$, if $K$ is \ec\ in $L$ as a field, then $(K,d)$ is \ec\ in $(L,\delta)$. \end{FACT} \noindent Examples of differentially large fields in characteristic $p>0$ are differentially closed fields in the sense of Wood \cite{Wood1973}, and also separably differentially closed fields in the sense of Ino and the first author \cite{InoLS2023}. Recall that a differential field $(K,\delta)$ is said to be separably differentially closed if for every differential field extension $(L,\delta)/(K,\delta)$ with $L/K$ separable (as fields), $(K,\delta)$ is \ec\ in $(L,\delta)$. To see that this class of differential fields is differentially large one only needs to note that if $K$ is \ec\ in $L$ as a field, then $L/K$ is separable. \par \medskip Let $(K,\delta)$ be a differential field (of arbitrary characteristic). In what follows we freely and interchangeably view any differential polynomial $f\in K\{x\}$ of order $n$ as a differential polynomial in the differential variables $x=(x_1,\dots,x_m)$ and also as a polynomial in $m(n+1)$ algebraic variables $x,\delta x,\dots,\delta^n x$. It will be clear from the context which view we are taking; for instance, if $a\in K^{m(n+1)}$ and we write $f(a)=0$, we mean viewing $f$ as a polynomial in $m(n+1)$ variables. \par \medskip In Theorem~\ref{DLAxiomsOrdinary} below we provide Blum-style axioms for ordinary differentially large fields of arbitrary characteristic. The proof relies on the following fact and its consequences, about extending derivations. \begin{FACT}{Fact}{ExtendKaehlerDerivation} \cite[Theorem 18, \S IV.7]{Jacobs1964} Suppose $L/K$ is a separable field extension. If $\delta:K\to L$ is a derivation, then $\delta$ can be extended to a derivation $L\to L$. \end{FACT} \begin{FACT}{Corollary}{DiffToDiffAlg} Let $(K,\delta)\subseteq (L,\delta)$ be an extension of differential fields and let $E$ be a subset of $L$ with $L/K(E)$ separable. Then there is a derivation $\partial:K(E\cup\delta(E))\lra K(E\cup\delta(E))$ that restricts to $\delta$ on $K(E)$. \par If $E$ is finite, then for each such $\partial $ there is some $f\in K[E\cup\delta(E)]$ such that $\partial$ restricts to a derivation of the localisation $K[E\cup\delta(E)]_f\lra K[E\cup\delta(E)]_f$. \end{FACT} \begin{proof} Since $\delta(K(E))\subseteq K(E\cup\delta(E))$ we may apply \ref{ExtendKaehlerDerivation} to the derivation $\delta|_{K(E)}:K(E)\lra K(E\cup\delta(E))$ and get a derivation $\partial:K(E\cup\delta(E))\lra K(E\cup\delta(E))$ that restricts to $\delta$ on $K(E)$. Assume then that $E$ is finite. There is some nonzero $f\in K[E\cup\delta(E)]$ such that $f\mal \partial(\delta(a))\in K[E\cup\delta(E)]$ for each $a\in E$. Obviously $f$ has the required property. \end{proof} \begin{FACT}{Proposition}{ReduceToDiffAlgAlgebra} Let $K$ be a differential field and let $S=(S,\delta)$ be a differentially finitely generated $K$-algebra and a domain such that $S/K$ is separable, i.e., $\Quot (S)/K$ is a separable field extension. Let $A$ be a finitely generated $K$-subalgebra of $S$. Then there are an element $f\in S$, a finitely generated $K$-subalgebra $B$ of $S_f$ containing $A$, a derivation $\partial $ on $B$ and a differential $K$-algebra homomorphism $S\lra (B,\partial)$ that restricts to the identity map on $A$. In particular $\partial a=\delta a$ for all $a\in A$.\looseness=-1 \end{FACT} \begin{proof} Let $b\in S^n$ be such that $S$ is the differential $K$-algebra generated by $b$ and $A\subseteq K[b]$. Let $\Dp=\{f\in K\{x\}\st f(b)=0\}$ be the differential vanishing ideal of $b$ over $K$. Then $\Dp$ is a separable prime differential ideal; separability is due to fact that $K\{x\}/\Dp$ is $K$-isomorphic to $S$. By the differential basis theorem of Kolchin \cite[Corollary 4, \S III.5]{Kolchi1973}, there is a finite set $\Sigma\subseteq \Dp$ that generates $\Dp$ as a radical differential ideal. Take $d\geq 1$ such that each derivative of any $x_1,\ldots,x_n$ occurring in some polynomial from $\Sigma$ has order $\leq d$. Finally take $$E=\{\delta^kb_i\st i\in \{1,\ldots,n\},\ k\leq d\}\subseteq S.$$ By possibly taking a larger $d$, a result of Kolchin appearing in \cite[Lemma 1, \S III.2]{Kolchi1973} tells us that $S/K(E)$ is separable. By \ref{DiffToDiffAlg} there are $f\in K[E\cup\delta(E)]$ and a derivation $\partial$ of $B:=K[E\cup\delta(E)]_f$ that restricts to $\delta$ on $K[E]$. Then $\partial^kb_i=\delta^kb_i$ for all $i\in \{1,\ldots,n\},\ k\leq d$ and therefore $b$ is a solution to $\Sigma=0$ in $(B,\partial)$. Consequently, the identity map of $K\cup \{b_1,\ldots,b_n\}$ extends to a differential $K$-algebra homomorphism $\phi:S\lra (B,\partial)$. By choice of $b$, the map $\phi$ restricts to the identity map of $A$. \end{proof} \begin{FACT}{Corollary}{ReduceToDiffAlg} Let $\Sigma$ be a set of differential polynomials over $(K,\delta)$ in finitely many differential variables. Suppose $\Sigma=0$ has a solution in some differential field extension $(L,\delta)$ with $L/K$ separable. Then there is a finitely generated $K$-subalgebra $B$ of $L$ and a derivation $\partial$ of $B$ such that $(B,\partial)$ has a solution to $\Sigma=0$. In particular, $(B,\partial)$ is differentially algebraic over $(K,\delta)$ and $B/K$ is separable. \par Notice that if $K$ is \ec\ in $L$ as a field then $K$ is also \ec\ in $B$ as a field. \end{FACT} \begin{proof} By assumption, there is a solution of $\Sigma=0$ in a differentially finitely generated $K$-subalgebra $S$ of $L$. Now apply \ref{ReduceToDiffAlgAlgebra} to $S$ and $A=K$. \end{proof} \begin{FACT}{Remark}{} In the case of several commuting derivations statements similar to \ref{ReduceToDiffAlgAlgebra} and \ref{ReduceToDiffAlg} fail in general. This follows from examples produced by Johnson, Reinhart, and Rubel \cite[Theorem 2]{JoReRu1995}. In particular, working over $(\C(z_1,z_2),\delta_1\equiv\frac{\partial}{\partial z_1},\delta_2\equiv\frac{\partial}{\partial z_2})$, they prove that the PDE $$\delta_2(x)=\left(1- \frac{z_1}{z_2} \right)\, x +1$$ has no differential algebraic solutions (equivalently, has no solution in a differential field extension of finite transcendence degree over $\C$). \end{FACT} \smallskip\noindent Recall that for a differential polynomial $f\in K\{x\}$, where $x$ is a single differential variable, we denote by $s_f$ the separant of $f$; namely, the formal partial derivative of $f$ with respect to its highest order variable. Furthermore we write \[ [f]:s_f^\infty= \{g\in K\{x\}: s^m g\in [f] \; \text{ for some } m\geq 0\}. \] \begin{FACT}{Observation}{ReduceToIrred} Let $K$ be a differential field and let $f\in K\{x\}$ for $x$ a single differential variable. Let $n=\ord(f)\geq 0$ and let $a\in K^{n+1}$ with $f(a)=0$ and $s_f(a)\neq 0$. Then there is an irreducible factor $h$ of $f$ with $\ord(h)=n$, $h(a)=0$ and $s_h(a)\neq 0$.\looseness=-1 \end{FACT} \begin{proof} Let $f_0,f_1\in K\{x\}$, with $f_0$ irreducible, $f=f_0\mal f_1$ and $\ord(f_0)=n$. Then \[ (*)\qquad\qquad s_f=\frac{\partial f}{\partial x_n}=\frac{\partial f_0}{\partial x_n}\,\mal\, f_1+f_0\,\mal\, \frac{\partial f_1}{\partial x_n}.\qquad\qquad \] If $f_0(a)=0$, then $(*)$ implies $s_{f_0}(a)=\frac{\partial f_0}{\partial x_n}(a)\neq 0$. If $f_0(a)\neq 0$, then $f_1(a)=0$ and $(*)$ shows $s_{f_1}(a)=\frac{\partial f_1}{\partial x_n}(a)\neq 0$; hence also $\ord(f_1)=n$ and in this case we may replace $f$ by $f_1$ and proceed by induction. \end{proof} \par \noindent We now come to the promised axiomatisation. \par \goodbreak \begin{FACT}{Theorem}{DLAxiomsOrdinary} Let $(K,\delta)$ be an ordinary differential field of arbitrary characteristic. The following conditions are equivalent. \begin{enumerate}[(i),itemsep=1ex] \item $(K,\delta)$ is differentially large. \item\labelx{DLAxiomsOrdinaryOrder} $K$ is large as a field and for every pair $f,g\in K\{x\}$, where $x$ is a single differential variable, with $g$ nonzero and $\ord(f)>\ord(g)$, if the system \[ f(x)=0 \ \&\  s_f(x)\neq 0 \] \noindent has an algebraic solution in $K$, then $f(x)=0\ \&\ g(x)\neq 0$ has a differential solution in $K$.\footnote{By \ref{ReduceToIrred} we may also assume that $f$ is irreducible in this condition.} \item\labelx{DLAxiomsOrdinaryMany} For every pair $f,g\in K\{x\}$, where $x$ is a single differential variable, with $\ord (f)\geq 1$ and $\ord(f)\geq\ord(g)$, if the system \[ f(x)=0 \ \&\ g(x)\cdot s_f(x)\neq 0 \] has an algebraic solution in $K$, then it has infinitely many differential solutions in $K$. \par \end{enumerate} Notice that each of the properties \ref{DLAxiomsOrdinaryOrder} and \ref{DLAxiomsOrdinaryMany} gives an axiom scheme for a first order axiomatization of differential largeness in the language of differential rings. \end{FACT} \begin{proof} (i)$\Ra $\ref{DLAxiomsOrdinaryMany}. Let $f,g\in K\{x\}$ with $\ord (f)\geq 1$ and $\ord(f)\geq \ord(g)$ and assume \[ (\dagger)\qquad\qquad f(x)=0 \ \&\ g(x)\cdot s_f(x)\neq 0\qquad\qquad \] has an algebraic solution in $K$. Let $n=\ord(f)$. By \ref{ReduceToIrred}, we may assume that $f$ is irreducible. Let $\Dp=[f]:s_f^\infty$. Since $s_f\neq 0$, Theorem 3.1(2) of \cite{InoLS2023} says that $\Dp$ is a separable prime differential ideal of $K\{x\}$. We write $a=x\mod\, \Dp$. Now, an algebraic solution of $f(x)=0 \ \&\ s_f(x)\neq 0$ in $K$ is a smooth $K$-rational point of \[ K[x_0,\ldots,x_n]/(f)\cong_KK[a,\ldots,a^{(n)}]. \] The largeness of $K$ yields that $K$ is \ec\ in $K(a,\ldots,a^{(n)})$. Since the latter is equal to the differential field $K\langle a\rangle$ generated by $a$ over $K$, differential largeness implies that $(K,\delta)$ is \ec\ in $(K\langle a\rangle, \delta)$. \par Since $\ord(f)\geq \ord(g)$ and $(\dagger)$ has an algebraic solution in $K$, Lemma 3.6(1) of \cite{InoLS2023} implies that $g\cdot s_f\notin \Dp$. Hence $a$ is a differential solution of $(\dagger)$ in $K\<a\>$. As $(K,\delta)$ is \ec\ in $(K\langle a\rangle, \delta)$ also $K$ has a differential solution $\alpha$ of $(\dagger)$. To argue that there are infinitely many solutions, note that $g\cdot (x-\alpha)$ has again order at most $\ord(f)$. By largeness of $K$ and the assumption $\ord(f)\geq 1$, there is an algebraic solution of the new system where we replace $g$ with $g\cdot (x-\alpha)$. It follows, by repeating the above argument, that there are infinitely many differential solutions of $(\dagger)$ in $K$.\looseness=-1 \par \smallskip\noindent \ref{DLAxiomsOrdinaryMany}$\Ra $\ref{DLAxiomsOrdinaryOrder} It suffices to show that $K$ is large as a field. By \cite[Lemma 5.3.1, p. 67]{Jarden2011}, a field $K$ is large if and only if for every absolutely irreducible polynomial $F(X,Y)\in K[X,Y]$, if there is a point $(a,b)\in K^2$ with $F(a,b) = 0$ and $\frac{\partial F}{\partial Y}(a,b)\neq 0$, then there are infinitely many such points. \par So take an absolutely irreducible polynomial $F(X,Y)\in K[X,Y]$ and some $(a,b)\in K^2$ with $F(a,b) = 0$ and $\frac{\partial F}{\partial Y}(a,b)\neq 0$. Consider the differential polynomial $f(x)=F(x,x')$. Then $f(x)=0\ \&\ s_f(x)\neq 0$ has an algebraic solution in $K$, namely $(a,b)$. By \ref{DLAxiomsOrdinaryMany} there are infinitely many differential solutions in $K$. But then there are infinitely many solutions to $F(X,Y) = 0$ and $\frac{\partial F}{\partial Y}(X,Y)\neq 0$ in $K$ as well.\looseness=-1 \par \smallskip\noindent \ref{DLAxiomsOrdinaryOrder}$\Ra $(i). To prove differential largeness, let $F$ be a differential field extension of $K$ such that $K$ is \ec\ in $F$ as a field. Note that then $F/K$ is separable. We need to show that $K$ is \ec\ in $F$ as a differential field. Let $\Sigma$ be a system of differential polynomials in $n$ differential variables over $K$ and assume that $\Sigma=0$ has a solution $a\in F^n$. We may assume that $F=K\langle a\rangle$. By \ref{ReduceToDiffAlg} applied to $F$, we may assume that $F$ is differentially algebraic over $K$ (and $F/K$ remains separable). \par Condition \ref{DLAxiomsOrdinaryOrder} guarantees that $[K:C_K]$ is infinite; hence, by the differential primitive element theorem \cite[Proposition 9, \S II.8, p.103]{Kolchi1973}, the differential field $F$ is differentially generated over $K$ by a single element $b\in F$. Let $\Dp$ be the prime differential ideal of $K\{x\}$ associated to $b$. Note that $\Dp$ is separable (over $K$). \par Then, by Theorem 3.1(1) of \cite{InoLS2023}, $\Dp=[f]:s_f^\infty$ for $f\in \Dp$ irreducible of minimal rank. Write $a=(a_1,\dots,a_n)$ and let $f_i,g\in K\{x\}$ with $a_i=\frac{f_i(b)}{g(b)}$. By the differential division algorithm \cite[\S I.9]{Kolchi1973} there are $h\in K\{x\}$ reduced with respect to $f$ and some $r\geq 0$ with \[ (i_fs_f)^rg\equiv h\quad \mod[f]. \] Since $f(b)=0$ and $i_f(b)\mal s_f(b)\neq 0$ we get $i_f^r(b)s_f^r(b)g(b)=h(b)\neq 0$. Hence, we may replace $g$ by $h$ and $f_i$ by $(i_fs_f)^r\mal f_i$ if necessary and assume that $g$ is reduced with respect to $f$. Notice that $a_i\in K\{b\}_{g(b)}$. \par Now, since $K$ \ec\ in $F$ as a field, the system $f(x)=0 \ \&\ s_f(x)\neq 0$ has an algebraic solutions in $K$. By condition \ref{DLAxiomsOrdinaryOrder}, the set $$\{f=0\}\cup \{q\neq 0\st q\in K\{x\} \text{ is nonzero and } \ord (q)<\ord (f)\}$$ is finitely satisfiable in the differential field $K$. Hence there is an elementary extension $L$ of the differential field $K$ having a differential solution $c$ to $f(x)=0$ such that $q(c)\neq 0$ for all $q\in K\{x\}$ with $\ord (q)<\ord (f)$. Since $f$ is irreducible, it follows that $q(c)\neq 0$ for all $q\in K\{x\}$ that are reduced with respect to $f$. \par In particular $f(c)=0\ \&\ g(c)\neq 0$. Since $K\prec L$ there is some $d\in K$ with $f(d)=0\ \&\ g(d)\neq 0$. This means there is a differential $K$-homomorphism $(K\{x\}/\Dp)_{g\mod \Dp}\lra K$. By choice of $\Dp$ we have $(K\{x\}/\Dp)_{g\mod \Dp}\cong K\{b\}_{g(b)}$ as differential $K$-algebras. Since $K\{a_1,\ldots,a_n\}\subseteq K\{b\}_{g(b)}$ we obtain a differential $K$-algebra homomorphism $K\{a_1,\ldots,a_n\}\lra K$ and this corresponds to a differential solution of $\Sigma=0$ in $K^n$. \end{proof} \par \noindent When $K$ is real closed, the above theorem yields a new axiomatisation of the theory \CODF. A differential field $(K,\delta)$ is a model of \CODF\ if and only if it is an existentially closed model of the theory of ordered differential fields. Axioms for \CODF\ appear in \cite{Singer1978a}. While the axioms there make explicit reference to the order, our new axioms are purely in the differential field language, namely: \begin{FACT}{Corollary}{CODFnewAxioms} Let $(K,\delta)$ be a differential field. The following are equivalent. \begin{enumerate}[(i)] \item $(K,\delta)\models$ \CODF. \item $K$ is real closed and for every pair $f,g\in K\{x\}$, where $x$ is a single differential variable, with $g$ nonzero and $\ord(f)>\ord(g)$, if the system \[ f(x)=0 \ \& \ s_f(x)\neq 0 \] has an algebraic solution in $K$, then $f(x)=0$ $\&$ $g(x)\neq 0$ has a differential solution in $K$. \end{enumerate} \end{FACT} \medskip \par \noindent Notice that every field $K$ is algebraically closed in the large field $K((t))$, but not every field is large. In the differential phrasing this changes: \begin{FACT}{Corollary}{DiffAlgClos} Let $L/K$ be an extension of differential fields. If $K$ is differentially algebraically closed in $L$ and $L$ is differentially large, then $K$ is differentially large as well. \end{FACT} \begin{proof} We verify \ref{DLAxiomsOrdinary}\ref{DLAxiomsOrdinaryMany}. Take $f,g\in K\{x\}$, $x$ a single differential variable, with $\ord(f)\geq 1$ and $\ord(f)\geq \ord(g)$, and assume that $f(x)=0 \ \&\ g(x)\cdot s_f(x)\neq 0$ has an algebraic solution in $K$. Since $L$ is differentially large, it has infinitely many differential solutions to $f(x)=0 \ \&\ g(x)\cdot s_f(x)\neq 0$. But then each of these solutions is differentially algebraic over $K$. Hence all these solutions are in $K$. \end{proof} \begin{FACT}{Remark}{} We note that the condition of a differential field $(K,\delta)$ being differentially algebraically closed in some extension $(L,\delta)$ is quite strong. Arguably, being differentially algebraically closed in an extension is not quite the right differential analogue of being algebraically closed in the field sense. We do not know whether the assumption in \ref{DiffAlgClos} can be weakened to only assuming that $K$ is constrainedly closed in $L$ (namely, every finite tuple from $L$ which is constrained over $K$, in the sense of Kolchin \cite[\S III.10]{Kolchi1973}, is from $K$). \end{FACT} We conclude this section with a geometric characterisation of being differentially large. Namely, in terms of algebraic D-varieties. Recall that an algebraic D-variety over $K$ is a pair $(V,s)$ where $V$ is an algebraic variety over $K$ and $s:V\to \tau V$ is a section over $K$ of the prolongation of $V$ (see \cite[\S2]{KowPil2005}, for instance). The latter is the algebraic bundle $\pi:\tau V\to V$ with the characteristic property that for any differential field extension $(L,\delta)$ of $(K,\delta)$ we have that if $a\in V(L)$ then $(a,\delta a)\in \tau V$.\looseness=-1 \begin{FACT}{Corollary}{DiffLargeDvariety} Let $K$ be a large field of arbitrary characteristic and let $\delta$ be a derivation of $K$. The following conditions are equivalent. \begin{enumerate}[(i)] \item $(K,\delta)$ is differentially large \item Let $V$ and $W$ be $K$-irreducible algebraic varieties with $W\subseteq \tau V$. If $\pi|_W:W\to V$ is a separable morphism and $W$ has a smooth $K$-point, then the set $$\{(a,\delta a)\in W:\; a\in V(K)\}$$ is Zariski dense in $W$. \item\labelx{DiffLargeDvarietyZarDense} Let $(V,s)$ be a $K$-irreducible algebraic D-variety. If $V$ has a smooth $K$-point, then the set $$\{a\in V(K): \; s(a)=(a,\delta(a))\}$$ is Zariski dense in $V$. \end{enumerate} \end{FACT} \begin{proof} (i)$\Rightarrow$(ii) Let $(a,b)$ be a $K$-generic point of $W$. Since $\pi_W:W\to V$ is a separable morphism, we obtain that $a$ is $K$-generic in $V$ and $K(a,b)/K(a)$ is a separable extension. Since $W\subseteq \tau V$, there is a derivation $\delta:K(a)\to K(a,b)$ extending the one on $K$ such that $\delta(a)=b$. As $K(a,b)/K(a)$ is separable, by \ref{ExtendKaehlerDerivation}, we can extend the derivation to $K(a,b)\to K(a,b)$. Then, for any nonempty Zariski-open $O_W\subseteq W$ over $K$, in the differential field extension $(K(a,b),\delta)$ we can find a solution to $x\in V$ and $(x,\delta x)\in O_W$ (namely, the tuple $a$). Since $W$ has a smooth $K$-point, we get that $K$ is \ec\ in $K(W)=K(a,b)$ as a field. By differential largeness, $(K,\delta)$ is \ec\ in $(K(a,b),\delta)$, and so we can find the desired solution in $K$. \par \medskip \par (ii)$\Rightarrow$\ref{DiffLargeDvarietyZarDense} If we let $W=s(V)\subseteq \tau V$, then the pair $V$ and $W$ satisfy the conditions of (ii) (note that if $b$ is a smooth point of $V$ then $(b,s(b))$ is a smooth point of $W$). If follows that the set of points in $W$ of the form $(a,\delta a)$ with $a\in V(K)$ is Zariski dense in $W$. But then, as $W=s(V)$, the set of points $a\in V$ such that $s(a)=(a,\delta a)$ must be Zariski dense in $V$. \par \medskip \par \ref{DiffLargeDvarietyZarDense}$\Rightarrow$(i) We verify \ref{DLAxiomsOrdinary}\ref{DLAxiomsOrdinaryOrder}. Let $f,g\in K\{x\}$ with $\ord(g)<\ord(f)$ and $g$ nonzero. Assume the system \[ f(x)=0 \ \&\  s_f(x)\neq 0 \] \noindent has an algebraic solution in $K$. In particular, $s_f\neq 0$. By Observation~\ref{ReduceToIrred}, we may assume that $f$ is irreducible. By Theorem 3.1(1) of \cite{InoLS2023}, $\Dp=[f]:s_f^\infty$ is a separable prime differential ideal of $K\{x\}$. Let $a=x+\Dp$ in the fraction field of $K\{x\}/\Dp$. Letting $n=\ord(f)$, we see that $(a,\delta a,\dots,\delta^{n-1}a)$ is algebraically independent over $K$ and $\delta^n a$ is separably algebraic over $K(a,\dots,\delta^{n-1}a)$. It follows that $$\delta^{n+1}a=\frac{h(a,\delta a, \dots, \delta^n a)}{s_f(a)}$$ for some $h\in K[t_0,\dots,t_n]$. Let $V$ be the localisation at $g\cdot s_f$ of the Zariski-locus of $(a,\delta a, \dots,\delta^n a)$ over $K$. From the assumptions (on existence of an algebraic solution in $K$), we see that $V$ has a smooth $K$-rational point and that the morphism on $V$ induced by $$(t_0,t_1, \dots,t_n)\mapsto ((t_0,t_1, \dots,t_n),(t_1,t_2,\dots,t_n, \frac{h(t_0,t_1, \dots, t_n)}{s_f(t_0,t_1, \dots, t_n)})$$ yields a regular algebraic map $s:V\to \tau V$. This equips $V$ with a D-variety structure. Then, the assumption of (iii) yields $\alpha\in V(K)$ such that $s(\alpha)=(\alpha,\delta \alpha)$. But then $\alpha$ is the desired differential solution of $f(x)=0\ \&\ g(x)\neq 0$ in $K$. \end{proof} \section{Power series in characteristic zero}\label{powerseries} \noindent In this section we assume fields are of characteristic zero, and thus the results on differentially large fields from \cite{LeSTre2024} may be deployed. We prove, in \ref{ECinDiffAlgSeries}, two further characterisations of being differentially large. \par For a differential field $K$ (ordinary throughout) we endow $K((t))$ with its natural derivation extending the given derivation on $K$ and satisfying $\delta(t)=1$; that is, \[ \delta(\sum_{n\geq k}a_nt^n)=\sum_{n\geq k}\delta(a_n)t^n+\sum_{n\geq k}na_nt^{n-1}. \] \par \noindent In \LST{CharDlargeI} it is shown that $(K,\delta)$ is differentially large if and only if $(K,\delta)$ is \ec\ in $(K((t)),\delta)$. We do not know if this characterisation extends to positive characteristic, the proof relies on the existence of a \textit{twisted version} of the Taylor morphism \LST{twistedTaylor}, whose construction picks up rational denominators. Below we prove that it suffices to ask for $(K,\delta)$ to be \ec\ in the differential subfield of $(K((t),\delta)$ consisting of differential algebraic elements (over $K$). \begin{FACT}{Definition}{} Let $K$ be a differential field and let $S$ be a differential $K$-algebra. We write $S_\mathrm{diffalg}$ for the differential subring of all $a\in S$ that are differentially algebraic over $K$. \end{FACT} \begin{FACT}{Remark}{} Since $K((t))$ is the localization of $K[[t]]$ at $t$, the fraction field of $K[[t]]_\mathrm{diffalg}$ is $K((t))_\mathrm{diffalg}$. \end{FACT} \begin{FACT}{Proposition}{DiffalgPowerSeries} Let $(K, \delta)$ be a differential field (of characteristic zero) that is large as a field and let $S$ be a differentially finitely generated $K$-algebra. If there is a $K$-algebra homomorphism $S\to L$ for some field extension $L/K$ in which $K$ is \ec\ (as a field, there are no derivations on $L$ given), then there is a differential $K$-algebra homomorphism $S\to K[[t]]_\mathrm{diffalg}$. \end{FACT} \begin{proof} By \LST{DiffPointFromPointInECext} there is a differential $K$-algebra homomorphism $\psi:S\to K[[t]]$. Applying \ref{ReduceToDiffAlgAlgebra} to $\psi(S)$ we may then find a finitely generated $K$-subalgebra $B$ of $K((t))$, a derivation $\partial $ of $B$ extending $\delta $ on $K$ together with a differential $K$-algebra homomorphism $\phi:\psi(S)\lra (B,\partial)$. By \LST{DiffPointFromPointInECext} applied to $(B,\partial)$ and the inclusion map $B\into K((t))$ there is a differential $K$-algebra homomorphism $\gamma:B\to K[[t]]$. Since $B$ is a finitely generated $K$-algebra, the image of $\gamma$ is in $K[[t]]_\mathrm{diffalg}$. Hence the map $\gamma\circ\phi\circ\psi:S\lra K[[t]]_\mathrm{diffalg}$ has the required property.\looseness=-1 \end{proof} \par \smallskip\noindent A special case of \ref{DiffalgPowerSeries} resembles an approximation statement over large and differential fields in the spirit of \cite[Theorem 2.1]{DenLip1984}: \begin{FACT}{Corollary}{WeakApproximation} Let $(K,\delta)$ be a differential field of characteristic zero such that $K$ is large as a field. Let $\Sigma$ be a system of differential polynomials in finitely many differential variables over $K$. If the differential ideal generated by $\Sigma$ has an algebraic solution in $K((t))$, then $\Sigma=0$ has a differential solution in $K[[t]]_\mathrm{diffalg}$. \end{FACT} \begin{proof} Apply \ref{DiffalgPowerSeries} to the differential coordinate ring of $\Sigma$. \end{proof} \begin{FACT}{Corollary}{ECinDiffAlgSeries} Let $K$ be a large field of characteristic 0 and let $\delta$ be a derivation of $K$. The following conditions are equivalent. \begin{enumerate}[(i)] \item $(K,\delta)$ is differentially large. \item $K$ is \ec\ in $K[[t]]_\mathrm{diffalg}$ as a differential field. \item For every $K$-irreducible algebraic D-variety $(V,s)$, if $V$ has a $K$-point, then there is $a\in V(K)$ such that $s(a)=(a,\delta a)$. \end{enumerate} \end{FACT} \begin{proof} (i)$\Ra $(ii) is a consequence of \LST{CharDlargeI}{IexClosedInSeries}, which says that $K$ is \ec\ in $K((t))$ as a differential field. \par \smallskip\noindent (ii)$\Ra $(i). By \ref{WeakApproximation} one verifies that $K$ is \ec\ in $K((t))$ as a differential field. Hence by \LST{CharDlargeI}, $(K,\delta)$ is differentially large \par \smallskip\noindent (iii)$\Ra $(i) We verify \ref{DiffLargeDvariety}\ref{DiffLargeDvarietyZarDense}. Let $(V,s)$ be a $K$-irreducible $D$-variety with a smooth $K$-point. Let $h\in K[V]$ nonzero. Then, there is an induced D-variety structure in the localisation $K[V]_h$. Denote this D-variety by $(W,t)$. As $K$ is large and $V$ has a smooth $K$-point, we get that $K$ is Zariski dense in $V$. Thus, $W$ has a $K$-point. The assumption now yields a $K$-point $b$ in $W$ such that $s(b)=(b,\delta b)$. As $h$ was arbitrary, it follows that the set of points $\{a\in V(K): s(a)=(a,\delta a)\}$ is Zariski dense.\looseness=-1 \par \smallskip\noindent (i)$\Ra $(iii) Let $(V,s)$ be a $K$-irreducible D-variety with a $K$-point. Applying \ref{DiffalgPowerSeries} with $S=K[V]$ and $L=K$, we find a $K((t))$-rational point $b$ of $V$ such that $s(b)=(b,\delta b)$. As $K$ is differentially large, it is \ec\ in $K((t))$ as a differential field. Hence, we can find such a point in $K$. \end{proof} \par \medskip\noindent We may now improve the construction of differentially large fields from \LST{ConstructionPowerSeries} in the ordinary case. A few preparations are necessary. \begin{FACT}{Proposition}{ConstructionDirectLimitOrdinary} Let $(K_i,f_{ij})_{i,j\in I}$ be a directed system of differential fields and differential embeddings with the following properties. \begin{enumerate}[(a)] \item\labelx{ConPropLargeOrdinary} All $K_i$ are large as fields. \item\labelx{ConPropECOrdinary} All embeddings $f_{ij}:K_i\lra K_j$ are isomorphisms onto a subfield of $K_j$ that is \ec\ in $K_j$ as a field. \item\labelx{ConPropSeriesOrdinary} For all $i\in I$ there exist $j\geq i$ and a differential homomorphism $K_i[[t]]_\mathrm{diffalg}\lra K_j$ extending $f_{ij}$. \end{enumerate} Then the direct limit $L$ of the directed system is a differentially large field. \end{FACT} \begin{proof} The proof is identical to the proof of \LST{ConstructionDirectLimit}, except we use \ref{DiffalgPowerSeries} in that proof instead of \LST{DiffPointFromPointInECext}. \end{proof} \begin{FACT}{Observation}{DiffAlgLocallyHenselian} Let $K$ be a differential field. Then $K[[t]]_\mathrm{diffalg}$ is a Henselian valuation ring. \end{FACT} \begin{proof} We write $S=K[[t]]_\mathrm{diffalg}$. Then $S$ is a valuation ring because if $f\in K((t))_\mathrm{diffalg}$ then the degree of $f$ is $\geq 0$, hence $f\in S$, or the degree of $f$ is negative and then $f^{-1}\in S$. Clearly the maximal ideal of $S$ is $t\mal S$. To verify that $S$ is Henselian it suffices to show that for all $\mu_2,\ldots,\mu_n\in \Dm$ there is some $f\in S$ with \[ 1+f+\mu_2 f^2+\ldots+\mu_n f^n=0. \] As $K[[t]]$ is Henselian, there is such an $f$ in $K[[t]]$. Obviously, $f\in S$. \end{proof} \begin{FACT}{Theorem}{ConstructionDiffAlgPowerSeries} Let $(K,\delta)$ be any differential field of characteristic zero. Set $K_0=K$ and let $K_{n+1}=K_n((t_{n}))_\mathrm{diffalg}.$ Then \( \bigcup_{n\geq 0} K_n \) is differentially large. \end{FACT} \begin{proof} By \ref{DiffAlgLocallyHenselian}, $K_n[[t_n]]_\mathrm{diffalg}$ is a Henselian valuation ring. By \cite{Pop2010}, $K_n((t_n))_\mathrm{diffalg}$ is a large field. We see that all assumptions of \ref{ConstructionDirectLimitOrdinary} are satisfied for the $K_n$ and the inclusion maps $K_n\into K_{n+k}$. Now the argument for \LST{ConstructionPowerSeries}{ConstructionPowerSeriesA} can be copied, where we use \ref{ConstructionDirectLimitOrdinary} instead of \LST{ConstructionDirectLimit}. \end{proof} \begin{FACT}{Remark}{} When $K$ is furnished with the trivial derivation, then the derivation on $K((t))$ is the standard one, namely $\frac{\dd}{\dd t}$, and hence by standard results on differential transcendence we know that $K((t))_\mathrm{diffalg}$ is a \textit{proper} differential subfield of $K((t))$. For instance, by H\"older's theorem, for $K=\C$, the Gamma function is in the latter while not in the former. We do not know whether this proper containment holds for an arbitrary differential field $K$. \end{FACT} We conclude this section by discussing possible improvements of \ref{DiffalgPowerSeries}. \begin{FACT}{Counterexample}{} If $K$ is algebraically closed in \ref{DiffalgPowerSeries} then a stronger conclusion holds, namely there is a differential $K$-algebra homomorphism $S\to K[[t]]$ whose image is constrained. The reason is that there is a differential homomorphism $\eps:S\lra K^\mathrm{diff}$ and then one can apply \ref{DiffalgPowerSeries} to obtain a differential embedding of the image of $\eps$ into $K[[t]]$. \par However, if $K$ is not algebraically closed then in general there is no differential $K$-algebra homomorphism $S\to K[[t]]$ whose image is constrained. To see an example, consider the ordered field $\R(z)$ where $z>\R$ and let $K$ be its real closure. We furnish $K$ with the unique derivation extending the standard derivation $\frac{\dd}{\dd z}$ on $\R(z)$. Let $x$ be a new transcendental element and let $R$ be the real closure of the ordered field $K(x)$ with the ordering $x>K$. Extend the derivation of $K$ to $R$ by setting $\delta(x)=0$. Let $y$ be a square root of $x-z$ in $R$ and let $S$ be the differential $K$-subalgebra of $R$ generated by $y$, hence $S=K[y,y^{-1}]$. Now if $\phi:S\lra K[[t]]$ is a differential $K$-algebra homomorphism, then $\phi(x)'=\phi(x')=0$ and $\phi(x)=\phi(y^2+z)=\phi(y)^2+z$, hence $\phi(x)$ is a constant and $\phi(x)-z$ is a square. As $z>\R$ and $\R$ is the constant field of $K$, we see that $\phi(x)$ cannot be in $K$. Hence $\phi(x)$ is a new constant of $K[[t]]$ and therefore it is not constrained over $K$. \end{FACT} \begin{FACT}{:On the canonicity of differentially algebraic solutions}{} Let $K$ be a differential field. If $S$ is a differentially finitely generated $K$-algebra and $\phi:S\lra K$ is a $K$-algebra homomorphism, then by \LST{DiffPointFromPointInECext} one can explicitly construct a differential $K$-algebra homomorphism $\psi:S\lra K[[t]]$, namely one can take $\psi$ to be the twisted Taylor morphism $T^*_\phi$ associated to $\phi$. Now, by \ref{DiffalgPowerSeries} there is even a differential $K$-algebra homomorphism $\rho:S\lra K[[t]]_\mathrm{diffalg}$ and one might ask whether $\rho$ can also be obtained in some canonical form out of $\phi$. However Gabriel Ng has shown that this is not possible. We refer to \cite[Proposition 7.11]{Ng2023} for details. \end{FACT} \section{Expansions of large fields to a differentially large field}\label{expansionslarge} \noindent The main goal of this section is \ref{DLexistence} which implies that any large field of characteristic zero of infinite transcendence degree over $\Q$ can be expanded to a differentially large field. A further consequence of \ref{DLexistence} is \ref{CODFprimModels}, which says that prime model extensions in \CODF\ only exist in the trivial case. Throughout this section fields are assumed to be of characteristic zero.  \begin{FACT}{Notation}{DefnDLproblem} Let $K$ be a field (of characteristic zero). A \notion{differentially large problem} of $K$ is a pair $(f,g)$ of polynomials from $K\{x\}=K[x_0,x_1,\dots]$ such that $f$ is of order $n\geq 0$, the order of $g$ is strictly less than $n$ and for which there is an element $(c_0,\ldots,c_n)\in K^{n+1}$ such that \[ f(c_0,\ldots,c_n)=0\ \&\ s_f(c_0,\ldots,c_n)\neq 0. \] We call ${\bar c}$ an \notion[]{algebraic solution of the differentially large problem}. Obviously a differentially large problem over $K$ remains a differentially large problem over every field extension of $K$. If $d$ is a derivation of $K$, then a \notion{solution of a differentially large problem} of $K$ in a differential field $(L,\delta)$ extending $(K,d)$ is an element $a\in L$ with $f(a)=0\ \&\ g(a)\neq 0$, where polynomials are now evaluated as differential polynomials. \par \end{FACT} \begin{FACT}{Proposition}{SolveDLProblem} Let $L/K$ be a field extension, $n\in\N$ and assume that $\trdeg(L/K)\geq n$. Let $(f,g)$ be a differentially large problem of $K$ with $\ord(f)=n$. Let $d$ be a derivation of $K$ and assume $L$ is large. \par Then there is a subfield $K_1$ of $L$ that is finitely generated over $K$ as a field, a derivation $\delta$ of $K_1$ extending $d$ and a solution $a\in K_1$ of the differentially large problem $(f,g)$ such that $a,\delta a,\ldots,\delta^{n-1}a$ are algebraically independent over $K$. \end{FACT} \begin{proof} Let ${\bar x}=(x_0,\ldots,x_n)$ and let $Z$ be the solution set in $L$ of the system \[ f({\bar x})=0\ \&\ s_f({\bar x})\neq 0\ \&\ g({\bar x})\neq 0. \] \claim There exists a point $(a_0,\ldots,a_n)\in Z$ with $\trdeg(a_0,\ldots,a_{n}/K)=n$. \par \noindent \textit{Proof.} Let $W$ be the variety defined by the two polynomials \[ f({\bar x}),\ y\mal s_f({\bar x})\mal g({\bar x})-1\in K[{\bar x},y]. \] Write $h({\bar x},y)=y\mal s_f ({\bar x})\mal g({\bar x})-1$. Then any common zero $({\bar a},c)$ of $f$ and $h$ in the algebraic closure of $L$ is a regular point of $W$, because $c\mal s_f ({\bar a})\mal g({\bar a})-1=0$ implies $\frac{\partial f}{\partial x_n}({\bar a})\neq 0$ and obviously $\frac{\partial h}{\partial y}=s_f \mal g$ does not vanish at ${\bar a}$. Hence the determinant of the matrix $\begin{pmatrix} \frac{\partial f}{\partial x_n} & \frac{\partial f}{\partial y} \\ \frac{\partial h}{\partial x_n} & \frac{\partial h}{\partial y} \end{pmatrix}$ is not zero at $({\bar a},c)$. This shows that $W$ is smooth. \par Since $(f,g)$ is a differentially large problem of $K$ we know that $W$ has a $K$-rational point. By \cite[Theorem 1]{Fehm2011}, using $\trdeg(L/K)\geq n=\dim(W)$, there is a $K$-embedding $K(W)\lra L$. A generic point of $W$ in $K(W)$ is then mapped to a point $(a_0,\ldots,a_n)\in Z$ with $\trdeg(a_0,\ldots,a_{n}/K)=n$. \hfill$\diamond$ \par As $s_f (a_0,\ldots,a_n)\neq 0$, $a_n$ is algebraic over $K(a_0,\ldots,a_{n-1})$. But now we see that $K_1:=K(a_0,\ldots,a_n)$ is isomorphic to the quotient field of $K\{x\}/\Dp$, where $\Dp=[f]:s_f^\infty$. This induces a derivation $\delta$ on $K_1$ and this derivation has the required properties: $a=a_0$ solves the given differentially large problem. \end{proof} \begin{FACT}{Theorem}{DLexistence} Let $L/K$ be an extension of fields of characteristic $0$ and suppose $L$ is a large field. Let $d$ be a derivation of $K$. If $\trdeg(L/K)\geq \card(K)$, then there is a derivation $\delta$ of $L$ extending $d$ such that $(L,\delta)$ is differentially large. \par [Under necessary assumptions on the constant field $C$ of $K$ we will show in Theorem \ref{ExtendToDL} that we may in addition find such a $\delta$ whose constant field is algebraic over $C$.] \end{FACT} \begin{proof} Let $\kappa=\card(K)$. By extending $K$ and $d$ we may assume that $\trdeg(L/K)=\kappa$. Let $\{t_i\st i<\kappa\}$ be a transcendence basis of $L$ over $K$ and let $(f_i,g_i)_{i\in \kappa }$ be a list of all differentially large problems of $L$; so here $f_i,g_i\in L\{x\}$ in the terminology of \ref{DefnDLproblem}. \par For $i<\kappa$ we define a subfield $K_i$ of $L$ and a derivation $d_i$ of $K_i$ such that \begin{enumerate}[(a)] \item $K_i$ contains $t_i$, $\trdeg(K_i/K)$ is finite for finite $i$ and $\trdeg(K_i/K)\leq \card(i)$ for $i\geq \omega$, \item $(K_i,d_i)$ extends $(K_j,d_{j})$ for $j<i$, and \item $(K_i,d_i)$ solves the differentially large problem $(f_i,g_i)$. \end{enumerate} \par \noindent Suppose $i<\kappa$ and $(K_j,d_j)$ has already been defined with properties (a)--(c); this also covers the case $i=0$. Let ${\bar b}\subseteq L$ be finite with $f_i,g_i\in K(\bar b)\{x\}$ such that there is an algebraic solution of the differentially large problem $(f_i,g_i)$ in $K({\bar b})$. Let $K_*$ be the field generated by $K(t_i,{\bar b})\cup \bigcup_{j<i}K_j$ and extend the derivation $\bigcup_{j<i}d_j$ to a derivation $d_*$ of $K_*$ arbitrarily. Obviously then $\trdeg(K_*/K)$ is finite if $i$ is finite and $\leq \card(i)$ otherwise. \par Consequently $\trdeg(L/K_*)$ is infinite and we may apply \ref{SolveDLProblem} to the extension $K_*\subseteq L$, the derivation $d_*$ and the differentially large problem $(f_i,d_i)$. We obtain an extension $(K_{i},d_{i})$ of $(K_*,d_*)$ such that $K_{i}$ is a subfield of $L$ that is finitely generated over $K_*$. Clearly $(K_{i},d_{i})$ satisfies (a)--(c). \par \smallskip\noindent Then $L=\bigcup_{i<\kappa}K_i$ and by \ref{DLAxiomsOrdinary} the differential field $(L,\partial)$ with $\partial=\bigcup_{i<\kappa}d_i$ is differentially large. \end{proof} \begin{FACT}{Remark}{} In characteristic $p>0$ the conclusion in \ref{DLexistence} fails even under the assumption that $L/K$ is separable. For example $L$ might be perfect (as a field), and hence any derivation on $L$ is trivial. \end{FACT} \begin{FACT}{Corollary}{AllLargeSupportDL} A large field $L$ of characteristic zero is of infinite transcendence degree if and only if there is a derivation $d$ of $L$ such that $(L,d)$ is differentially large.\looseness=-1 \end{FACT} \begin{proof} If $L$ has infinite transcendence degree, then by \ref{DLexistence} applied with $K=\Q$ shows that there is a derivation $d$ of $L$ such that $(L,d)$ is differentially large. For the converse assume there is a derivation $d$ of $L$ such that $(L,d)$ is differentially large. By \LST{algextdifflarge}, the algebraic closure $\overline{L}$ of $L$ is a DCF. We may then replace $L$ by the differential closure of $\Q$. By the non-minimality of the differential closure of $\Q$ (\cite{Rosenl1974}), there is an embedding $L\lra L$ that is not surjective. Hence $L$ cannot have finite transcendence degree. \end{proof} \section{Differentially large fields with prescribed constant field} \label{SectionExpansionsWithConstants} \par \noindent We now aim to provide a version of Theorem~\ref{DLexistence} for all differentially large fields of characteristic 0 without extending the constants. More precisely, we prove in Theorem~\ref{ExtendToDL} below that for a field extension $L/K$ with $\trdeg(L/K)\geq \card(K)$, if $d$ be a derivation of $K$ whose constant field $C_K$ is dense in $L$ for the \'etale open topology of $L$, then there is an extension $\delta$ of $d$ on $L$ such that $(L,\delta)$ is differentially large and $C_L$ is algebraic over $C_K$. Hence, under a density assumption of the constants $C_K$, if $C_K$ is algebraically closed in $L$, then the construction of the derivation in Theorem~\ref{DLexistence} can be performed without introducing new constants. We conclude this section with an application to dense pairs of large fields in \ref{ExamplesCompleteTheoriesOfPairs}. \par \smallskip\noindent We first (briefly) introduce the notion of \textit{$L$-prime ideals} and tuples in the context of a fixed field extension $L/K$. For a differential ring $S$ we write $\Sped(S)$ for the subspace of $\Spec(S)$ consisting of the differential prime ideals of $S$. \begin{FACT}{Definition}{DefnLconstrained} Let $L/K$ be an extension of fields and let $d$ be a derivation of $K$. Let $S$ be a differential $K$-algebra. We call a prime ideal $\Dp$ of $S$ a \notion{differential $L$-prime ideal} if it is differential and $S/\Dp$ can be embedded into $L$ as a $K$-algebra; observe that there is no derivation given on $L$. We write $\Sped_L(S)$ for the subspace of $\Sped(A)$ consisting of differential $L$-prime ideals. We say that a point $\Dp\in \Sped(S)$ is \notion{$L$-locally closed} if it is a locally closed point of $\Sped_L(S)$. If $S$ is finitely generated as a $K$-algebra and $(0)$ is the unique $L$-locally closed point of $\Sped(S)$ we say that $S$ is \notion{$L$-simple}. Note that $L$-simplicity implies that there is a $K$-algebra embedding $S\lra L$.\looseness=-1 \end{FACT} \begin{FACT}{Examples}{} \begin{enumerate}[(i)] \item If $L$ is an algebraically closed field of infinite transcendence degree and $S$ is a differentially finitely generated $K$-algebra, then $\Sped_L(S)=\Sped(S)$, and $L$-locally closed is the same as being \notion{constrained} in the sense of Kolchin \cite{Kolchi1974}.\looseness=-1 \item If $K$ is real closed and $L$ is an $|S|^+$-saturated real closed field, then $\Sped_L(S)$ is the subspace of differential prime ideals that are real $S$ (meaning: $-1$ is not a sum of squares of $S$). When we are in this example we will say \notion{real constrained} instead of $L$-locally closed. \item Clearly, being constrained and real implies real-constrained. However, the converse does not always hold. For instance, consider the real closure $K=\Q(t)^\mathrm{rcl}$, where $\Q(t)$ is equipped with the unique ordering such that $t>\Q$ and with the unique derivation extending $\frac{d}{dt}$ on $\Q(t)$. Let $\alpha_2$ be a transcendental over $K$. In the formally real field $K(\alpha_2)$ define a derivation $\delta$ that extends the one on $K$ such that $\delta(\alpha_2)=\frac{-1}{2\alpha_2}$. Let $\alpha_1=t+\alpha_2^2$. Then, $\delta(\alpha_1)=0$ and $\alpha_1$ is transcendental over $K$ (as $\alpha_1>\Q^\mathrm{rcl}$). Now consider $\alpha=(\alpha_1,\alpha_2)$. Clearly $\alpha$ is not constrained over $K$ (as $\alpha_1$ is a constant which is not algebraic over $C_K=\Q^\mathrm{rcl}$). But $\alpha$ is real-constrained over $K$. Indeed, for any differential specialisation $\beta=(\beta_1,\beta_2)$ of $\alpha$ over $K$ with $K\langle \beta\rangle$ formally real, we see that $\beta_1$ is transcendental over $K$ (for the same reason that $\alpha_1$ was); thus, the map $K(\alpha)\to K(\beta)$ fixing $K$ and mapping $(\alpha_1,\alpha_2)\mapsto (\beta_1,\beta_2)$ is a differential isomorphism (not necessarily preserving the orderings). \end{enumerate} \end{FACT} \begin{FACT}{Observation}{LsimpleLDL} Let $L$ be a differentially large field and let $K\subseteq L$ be a differential subfield. Then for every $L$-simple $K$-algebra there is a differential $K$-embedding $S\lra L$. \end{FACT} \begin{proof} As $S$ is $L$-simple, there is a $K$-algebra embedding $\phi:S\lra L$. Now $S\otimes_KL$ is a finitely generated differential $L$-algebra and $\phi$ extends to an $L$-algebra homomorphism $S\otimes_KL\lra L$. By differential largeness, \LST{CharDlargeI}{IPointImpliesDiffPoint} says that there is a differential $L$-algebra homomorphism $\psi:S\otimes_KL\lra L$. Composing $\psi$ with the natural map $S\lra S\otimes_KL$ gives a differential $K$-algebra homomorphism $S\lra L$. Since $S$ is $L$-simple, this map is an embedding. \end{proof} \begin{FACT}{Proposition}{DLisolationRealisedInAffine} Let $L/K$ be a field extension and let $d$ be a derivation of $K$. Let $S$ be a differentially finitely generated $K$-algebra. If $\Sped_L(S)$ is nonempty, then there are $L$-locally closed points $\Dp$ of $S$ and for each such point there is some $q\in S$ such that $(S/\Dp)_q$ is $L$-simple. \end{FACT} \begin{proof} The existence of an $L$-locally closed point follows from Noetherianity of $\Sped(S)$: The patch closure $X$ of $\Sped_L(S)$ is then also Noetherian. Any locally closed point of $X$ is then in $\Sped_L(A)$ as one checks without difficulty. \par Now take an $L$-locally closed points $\Dp$ of $S$. Hence there is some $q\in S$ such that the prime ideal $\Dp$ is maximal for inclusion in $D(q)\cap \Sped_L(S)$. In other words the zero ideal is the unique element of $\Sped_L(A)$, where $A=(S/\Dp)_q$; in particular there is a $K$-algebra embedding $A\lra L$. \par By \ref{ReduceToDiffAlgAlgebra} there are $a\in A$ and a (not necessary differential) finitely generated $K$-subalgebra $B$ of $A_a$ and a derivation $\partial$ of $B$ together with a differential $K$-algebra homomorphism $f:A\lra (B,\partial)$. In particular the kernel of $f$ is in $\Sped(A)$. Since $B$ is a $K$-subalgebra of $A_a$, the kernel of $f$ is even in $\Sped_L(S)$ and so the kernel of $f$ is $0$. Thus $f$ is an embedding and $\qf(A)$ is of finite transcendence degree over $K$. After localization of $A$ by some element of $A$ we even see that $A$ is finitely generated as a $K$-algebra. \end{proof} \par In what follows we will talk about the \textit{\'etale open topology} on $K$-rational points of $K$-varieties for a field $K$, cf. \cite[p. 4034]{JTWY2024}. Explicitly we will only need a few basic properties of the \'etale open topology of $K$ itself and we only record what we need later on: \begin{FACT}{:The \'etale open topology.}{DefnEtOpenOnK}Let $K$ be a field. \begin{enumerate}[(i)] \item We call a subset $U$ of $K$ \notion{standard \'etale open} if it is the image of the projection $K^2\lra K$ onto the first coordinate of a set of the form $\{(a,b)\in K^2\st P(a,b)=0\ \&\ Q(a,b)\neq 0\}$, where $P,Q\in K[x,y]$ such that $\frac{\partial}{\partial y}P$ is invertible in the localization of $K[x,y]/(P)$ at $Q$. In the terminology of \cite[Definition 3.5.38]{Poone2017} these sets are precisely the images of $K$-rational points of \notion{standard \'etale morphisms} defined over $K$ with codomain $\A^1$. The standard \'etale open sets form a basis of a topology on $K$ which is the \notion{\'etale open} topology, cf. \cite[Definition 3.5.38]{Poone2017}, \cite[p. 4037]{JTWY2024}. \item \cite{JTWY2024} The field $K$ is large if and only if the \'etale open topology is not discrete. If $K$ is algebraically closed, then the \'etale open topology is the Zariski topology. If $K$ is real closed, then the \'etale open topology is the order topology. If $K$ possesses a henselian valuation $v$, then the \'etale open topology is the open ball topology of $v$. \end{enumerate} \end{FACT} \begin{FACT}{Proposition}{LLocClConstantsUnchanged} Let $L/K$ be a field extension such that $L$ is large and let $d$ be a nontrivial derivation of $K$. Let $S$ be an $L$-simple differential $K$-algebra. If $K$ is algebraically closed in $L$ and $C_K$ is dense in $L$ for the \'etale open topology of $L$, then the constant field of $\qf(S)$ is $C_K$. \end{FACT} \begin{proof} By $L$-simplicity, we may assume that $S$ is a $K$-subalgebra of $L$. We write $F=\qf(S)$ and $\delta$ for the derivation of $F$. Since $S$ is $L$-simple, there is $g\in S$ such that the localisation $\Sped_L(S_g)$ only consist of the zero ideal. \par Now suppose $f=\frac{p}{q}\in F$ is a constant, thus $\delta(f)=0$. We aim to show that $f$ is in $C_K$. We work in the localisation $S_{g\cdot q}$, and we view it as the coordinate ring of an affine variety $V$ defined over $K$. Then $f$ yields an algebraic map $f:V(L)\to L$. The image $W=f(V(L))$\footnote{Formally: $W$ is the set of all $\eps(f)$, where $\eps:S_{g\cdot q}\lra L$ is a $K$-algebra homomorphism.} is $K$-definable in $L$ in the language of rings. \par \smallskip \case 1 $W$ is infinite. \par \smallskip\noindent Since $W$ is an existentially $L$-definable set it must have nonempty interior for the \'etale open topology of $L$ by \cite[Corollary A, p.613]{WalYe2023}. Since $C_K$ is dense in $L$, there is $\eps\in V(L)$ -- hence a $K$-algebra homomorphism $S_{g\cdot q}\lra L$ -- such that $c:=\eps(f)\in C_K$ and we claim that $f=c$. \par Since $\delta(f-c)=0$, the ideal $\Da:=(f-c)$ of $S_{g\cdot q}$ is differential and contained in the kernel $\Dq$ of $\eps$. Choose any extension of the derivation of $K$ to $L$ and let $T_\eps^*:S_{g\cdot q}\lra L[[t]]$ be the twisted Taylor morphism of $\eps$ and that derivation. We write $\Dq^\#$ for the kernel of $T_\eps^*$ and obtain a differential $K$-algebra embedding $S_{g\cdot q}/\Dq^\#\into L[[t]]$\footnote{We only need the derivation on $L$ to obtain some $K$-algebra embedding as asserted, the classical Taylor morphism would not deliver this.}. It follows that the $K$-variety $V_1$ defined by $S_{g\cdot q}/\Dq^\#$ has a smooth $L[[t]]$-rational point. Since $L$ is a large field, $V_1$ also has a smooth $L$-rational point. Since $S$ is a $K$-subalgebra of $L$ we see that $\trdeg(L/K)\geq \trdeg(S/K)\geq \trdeg((S_{g\cdot q}/\Dq^\#)/K)$. By \cite[Theorem 1]{Fehm2011} we know that there is a $K$-algebra embedding $S_{g\cdot q}/\Dq^\#\into L$. We have shown that $\Dq^\#$ is in $\Sped_L(A)$. But $(0)$ is the only differential $L$-prime ideal of $S_{g\cdot q}$, thus $\Dq^\#=(0)$. On the other hand, $\Da$ is a differential ideal and contained in the kernel $\Dq$ of $\eps$. This implies $\Da\subseteq \Dq^\#=(0)$, showing that $f-c=0$ as required. \par \smallskip \case 2 $W$ is finite. \par \smallskip\noindent Since $V$ is irreducible and defined over $K$ the assumption that $K$ is algebraically closed in $L$, implies that the variety $V\times_KL$ is also irreducible. Consequently $V(L)$ is an irreducible subset for the Zariski topology of the $L$-rational points of $V$. It follows that the image $W$ of $f$ seen as a map $V(L)\lra L$ is also irreducible. As $W$ is finite, $W$ is a singleton set. Hence $f$ is a constant algebraic map and thus $f=c$.\looseness=-1 \end{proof} \begin{FACT}{Corollary}{LLocClConstantsUnchangedCor} Let $L/K$ be and extension of differential fields and suppose $L$ is large as a field. Suppose $K$ is algebraically closed in $L$ and $C_K$ is dense in $L$ for the \'etale open topology of $L$. Let $a\in L^n$ and $g\in K\{x\}$ with $g(a)\neq 0$. \begin{enumerate}[(i)] \item There is an $L$-simple differential $K$-algebra $S$ such that the fraction field of $S$ has constant field $C_K$, together with a differential $K$-algebra homomorphism $\phi:K\{a\}\lra S$ satisfying $g(\phi(a))\neq 0$. \item If $L$ is differentially large, then there is some $b\in L^n$ with $g(b)\neq 0$ such that $C_{K\<b\>}=C_K$ together with a differential $K$-algebra homomorphism $K\{a\}\lra K\{b\}$ mapping $a$ to $b$. \end{enumerate} \par \end{FACT} \begin{proof} (i) Since $(0)$ is a differential $L$-prime ideal, there is a differential $L$-prime ideal $\Dp$ of $K\{a\}$ with $g\notin\Dp$ and some $s\in K\{a\}\setminus \Dp$ such that $S=(K\{a\}/\Dp)_{g(a)\mal s}$ is $L$-simple. In particular $S$ can be embedded as a $K$-algebra into $L$. By \ref{LLocClConstantsUnchanged}, the fraction field of $S$ has constant field $C_K$. \par (ii) follows from (i) and \ref{LsimpleLDL}. \end{proof} \begin{FACT}{Theorem}{ExtendToDL} Let $L/K$ be a field extension such that $L$ is large with $\trdeg(L/K)\geq \card(K)$. Let $d$ be a derivation of $K$ whose constant field $C$ is dense in $L$ for the \'etale open topology of $L$. Then there is an extension $\delta$ of $d$ on $L$ such that $(L,\delta)$ is differentially large whose constant field is algebraic over $C$. \end{FACT} \begin{proof}Using \ref{LLocClConstantsUnchangedCor} we adapt the strategy of the proof of \ref{DLexistence}. \par We may replace $K$ by its algebraic closure of $K$ in $L$, hence we need to find a derivation on $L$ extending $d$ with constant field $C$. Let $\kappa=\card(K)$. By extending $K$ and $d$ we may assume that $\trdeg(L/K)=\kappa$ and that $K$ is not a constant field. Let $(f_i,g_i)_{i\in \kappa }$ be a list of all differentially large problems of $L$, where $f_0=x,g_0=1$; so here $f_i,g_i\in L\{x\}$ in the terminology of \ref{DefnDLproblem}. \par For $i<\kappa$, starting with $K_0=K$ we define a subfield $K_i$ of $L$ that is algebraically closed in $L$ and a derivation $d_i$ of $K_i$ such that \begin{enumerate}[(a)] \item\labelx{ConstrExtendDL} $(K_i,d_i)$ extends $(K_j,d_{j})$ for $j<i$. \item\labelx{ConstrSizeDL} $\trdeg(K_i/K)\leq \max\{\aleph_0,\card(i)\}$. \item\labelx{ConstrInfiniteDL} $\trdeg(L/K_i)$ is infinite (which is implied by \ref{ConstrSizeDL} when $L$ is uncountable). \item\labelx{ConstrProblemDL} $f_i,g_i\in K_i\{x\}$ and $(K_i,d_i)$ solves the differentially large problem $(f_i,g_i)$. \item\labelx{ConstrConstDL} $(K_i,d_i)$ has constant field $C$. \end{enumerate} \par \noindent Suppose $0<i<\kappa$ and $(K_j,d_j)$ has already been defined for $j<i$ with properties (a)--\ref{ConstrConstDL}. Let $K_*=\bigcup_{j<i}K_j$ with derivation $d_*=\bigcup_{j<i}d_j$. Obviously then $\trdeg(K_*/K)\leq \max\{\aleph_0,\card(i)\}$. If $i$ is infinite, then $L$ is uncountable and of size $>\card(i)$, hence $\trdeg(L/K_*)$ is infinite. If $i$ is finite, then $K_*=K_{i-1}$ and $\trdeg(L/K_*)$ is infinite as well. \par Since $\trdeg(L/K_*)$ is infinite, there is a countable infinite set $T\subseteq L$ that is algebraically independent over $K_*$ such that $\trdeg(L/K_*(T))$ is infinite, such that $f_i,g_i\in K_*(T)^\alg\{x\}$ and such that there is an algebraic solution of the differentially large problem $(f_i,g_i)$ in $K_*(T)^\alg$ (here the superscript '$\alg$' stands for the algebraic closure in $L$). \par It follows that $\trdeg(L/K_*(T))\geq \card(K_*(T))$ in either case. Now the field $K_*(T)$ is isomorphic to the fraction field $K_*\<x\>$ of the differential polynomial ring $K_*\{x\}$ and therefore there is a derivation $\partial$ of $K_*(T)$ extending $d_*$ such that $(K_*(T),\partial)$ is $K_*$-isomorphic to $K_*\<x\>$ with its natural derivation. It follows that the constant field of $\partial$ is the constant field of $K_*$, which is $C$ by property \ref{ConstrConstDL} in the induction hypothesis. Hence we may extend $\partial$ to the algebraic closure $K_+$ of $K_*(T)$ in $L$ without extending the constants. \par Since $\trdeg(L/K_+)\geq \card(K_+)$ we may now apply \ref{DLexistence} and extend $\partial$ to a derivation on $L$ such that $(L,\partial)$ is differentially large. Since $C$ is dense in $L$ it is also dense in $K_+$. Since $(f_i,g_i)$ is a differentially large problem of $K_+$ by choice of $T$, we may apply \ref{LLocClConstantsUnchangedCor}(ii), which shows that $(f_i,g_i)$ has a differential solution $a$ in $(L,\partial)$ such that $(K_+\<a\>,\partial)$ has constant field $C$. We may then define $K_i$ to be the algebraic closure of $(K_+\<a\>,\partial)$ in $L$ and see that all conditions (a)--\ref{ConstrConstDL} are satisfied. \par Finally, $L=\bigcup_{i<\kappa}K_i$ because for each $b\in L$, the differentially large problem $(x-b,1)$ is solved in the union. Hence the \factname\ follows. \end{proof} \par \noindent From the description of the \'etale open topology in \ref{DefnEtOpenOnK} we see that every differentially large field has a dense constant field for that topology, hence the density assumption in \ref{ExtendToDL} is necessary. \begin{FACT}{Corollary}{ExtendToDLcor} \begin{enumerate}[(i)] \item For every real closed subfield $K$ of the field $\R$ with $\trdeg(\R/K)\geq 2^{\aleph_0}$ and any derivation $d$ of $K$ there is a derivation $\delta$ on $\R$ extending $d$ such that $(\R,\delta)$ is a \CODF\ whose constant field is the constant field of $(K,d)$. Recall that the \'etale open topology of $\R$ is the Euclidean topology of $\R$. \item For every p-adically closed subfield $K$ of the field $\Q_p$ with $\trdeg(\Q_p/K)\geq 2^{\aleph_0}$ and any derivation $d$ of $K$ there is a derivation $\delta$ on $\Q_p$ extending $d$ such that $(\Q_p,\delta)$ is differentially large, whose constant field is the constant field of $(K,d)$. Recall that the \'etale open topology of $\Q_p$ is the valuation topology of $\Q_p$ and $\Q$ is dense in $\Q_p$. \item\labelx{ExtendToDLcorCountable} For every pair $K\subseteq L$ of countable fields, if $K$ is large, $K$ is algebraically closed in $L$ and dense in $L$ for the \'etale open topology of $L$ with $\trdeg(L/K)$ infinite, then there is a derivation $\delta$ on $L$ such that $(L,\delta)$ is differentially large with constant field $K$. \end{enumerate} \end{FACT} \begin{FACT}{Example}{DensePowerSeries} Dense pairs of fields $C\subseteq L$ as required in \ref{ExtendToDL} -- hence $L$ is large, $C$ is dense for the \'etale open topology of $L$ with $\trdeg(L/C)\geq \card(C)$ -- also occur naturally in power series fields: Let $k$ be any field of characteristic 0. We work in the Henselian valued field $k((t^\Q))$ of generalized power series of $k$. Let $C$ be the algebraic closure of $k(t)$ in $k((t^\Q))$ and let $M$ be the completion of $C$ for the valued field $C$, hence $M$ is the subfield of $k((t^\Q))$ consisting of power series whose support is cofinal in $\Q$ and of order type $\omega$. \par Now for any $\Q$-linearly independent set $\Lambda \subseteq k$ of cardinality $\card(k)$ (a baby version of) Ax's solution to the functional Schanuel conjecture implies that the series $\exp(\lambda \mal t)$ with $\lambda\in\Lambda$ are algebraically independent over $k$. \par So if $L$ is a large subfield of $M$ (for example, $M$ itself), containing all the $\exp(\lambda \mal t)$, then the pair $L/C$ has the required properties. \end{FACT} \smallskip\noindent A further consequence of \ref{ExtendToDL} is an application to dense pairs of large fields: \begin{FACT}{Theorem}{ObtainCompleteTheoriesOfPairs} Let $T$ be a theory of large fields of characteristic 0 in the language $\SL$ of rings. Let $T_\mathrm{pair}$ be the $\SL(P)$ theory of proper pairs $K\subsetneq L$ of models of $T$ for which $K$ is algebraically closed in $L$ and $K$ is dense in the \'etale open topology of $L$; here $P$ is a new unary predicate.\footnote{Notice that all these conditions are first order in the language $\SL(P)$, even when the \'etale open topology is not a definable field topology. Density of the subfield is preserved by elementary equivalence, as one sees by using the description of the open sets in \ref{DefnEtOpenOnK}.} Let $T^\delta$ be the theory $T$ together with the theory of differentially large fields in the language $\SL(\delta)$ of ordinary differential rings. \par If $T^\delta$ is a complete theory, then also $T_\mathrm{pair}$ is a complete theory. \end{FACT} \begin{proof} If $(L,K)\models T_\mathrm{pair}$, then a standard compactness argument shows that there is an elementary extension $(L',K')\succ (L,K)$ such that $\trdeg(L'/K')$ is infinite. By Skolem-L\"owenheim downwards there is a countable elementary restriction $(L'',K'')\prec (L',K')$ such that $\trdeg(L''/K'')$ is infinite. \par Now take $(L_1,K_1),(L_2,K_2)\models T_\mathrm{pair}$. In order to show that the pairs $(L_1,K_1)$ and $(L_2,K_2)$ are elementary equivalent we may apply the argument above to each pair and assume that $L_i$ is countable and of infinite transcendence degree over $K_i$. By \ref{ExtendToDLcor} we may expand $L_i$ to a differentially large field $(L_i,\delta_i)$ with constant field $K_i$. Hence the completeness of $T^\delta$ implies that $(L_1,\delta_1)$ and $(L_2,\delta_2)$ are elementary equivalent as differential fields and therefore the pairs $(L_1,K_1)$ and $(L_2,K_2)$ are elementary equivalent as well. \end{proof} \begin{FACT}{Corollary}{ExamplesCompleteTheoriesOfPairs} Let $T$ be a complete and model complete theory of large fields of characteristic 0 in the language $\SL$ of rings, possibly extended by constants. Then the theory of proper dense pairs of models of $T$ is complete. \par This applies to topologically large fields in the sense of \cite{GuzPoi2010}, like for example real closed fields and p-adically closed fields, but also to complete theories of pseudo-finite fields (where we need the constants for model completeness). \end{FACT} \begin{proof} By \cite[Cor. 4.8(iii)]{LeSTre2024}, $T^\delta$ is model complete, hence \ref{ObtainCompleteTheoriesOfPairs} applies. \end{proof} \par \noindent For the case of real closed fields, \ref{ExamplesCompleteTheoriesOfPairs} is A. Robinson's theorem saying that the theory of proper pairs of dense real closed fields is complete. \par \section{Prime models in CODF are algebraic} \noindent We now apply \ref{DLexistence} to answer a question about prime model extensions for \CODF. Recall that a \CODF\ in the sense of Singer (cf. \cite{Singer1978a}) is the same as a differentially large field that is real closed as a pure field. In \cite{Singer1978a}, Singer shows that \CODF\ has no prime model, i.e. there is no \CODF\ that embeds into all other \CODF s.\footnote{Note that \CODF\ is model complete in the language of differential rings, i.e. every embedding of \CODF s is elementary.} We now show that in fact no differential and formally real field (i.e. it possesses an ordering) has a prime model extension for \CODF\footnote{A \textit{prime model extension of $K$ for \CODF} is a model $\hat K$ of \CODF\ having $K$ as a differential subfield such that $\hat K$ embeds over $K$ as a differential field into any other \CODF\ that has $K$ as a differential subfield.}, unless its real closure is already a \CODF. In particular, no formally real field equipped with the trivial derivation has a prime model extension in \CODF. The proof is essentially an application of \ref{DLexistence} together with the following purely field theoretic fact. \begin{FACT}{Proposition}{RCFnoEmbeding} Let $R$ be a real closed field and let $\kappa$ be its cardinality. Then, there are real closed fields $M,N$ containing $R$ of transcendence degree $\kappa $ over $R$ with the following property: If $S\supseteq R$ is a real closed field then $S$ can be embedded over $R$ into $M$ and into $N$ if and only if $\trdeg(S/R)\leq 1$ and $R$ is Dedekind complete in $S$.\looseness=-1 \end{FACT} \begin{proof} We take $M\supseteq R$ by successively adjoining infinitely large elements $a_\alpha$ for $\alpha<\kappa$. Hence $a_\alpha>R(a_\beta\st \beta<\alpha)$ in $M$ and $M$ is algebraic over $R(a_\alpha\st \alpha<\kappa)$. Then $R$ is Dedekind complete in $M$ and $M$ has transcendence degree $\kappa$ over $R$. \par \smallskip\noindent For $N$ we may take any real closed subfield of $R((t^\R))$ of transcendence degree $\kappa$ over $R$. Such fields exist because of the following reason: Let $\Lambda$ be a basis of the $\Q$-vector space $R$. Since $R$ is real closed, the cardinality of $\Lambda$ is $\kappa$. Then the set $\{\exp(\lambda \mal t)\st \lambda\in\Lambda\}\subseteq R[[t]]$ is an algebraically independent subset of $R((t^\R))$ over $R$: this is a baby case of Ax's positive solution to the functional Schanuel conjecture, but is not difficult to prove directly. Hence we may take $N$ as the real closure of $R(\exp(\lambda \mal t)\st \lambda\in\Lambda)$ in $R((t^\R))$. Clearly $N$ has transcendence degree $\kappa$ over $R$.\looseness=-1 \par \smallskip Since $R$ is Dedekind complete in $M$ and in $N$, any real closed field $S$ containing $R$ with $\trdeg(S/R)\leq 1$ in which $R$ is Dedekind complete, can be embedded into $M$ and into $N$. It remains to show that any real closed subfield $S$ of $M$ containing $R$ that can be embedded into $N$ over $R$ is of transcendence degree at most 1 over $R$; note that $R$ is Dedekind complete in $S$ because $R$ is Dedekind complete in $M$ (and in $N$). \par For a contradiction, suppose $S$ has transcendence degree $2$ over $R$. We furnish $M$ with the valuation whose valuation ring is the convex hull of $R$ in $M$. Real closures are now taken in $M$ throughout and this is indicated by the superscript~$^\mathrm{rcl}$. Take ${\bar a}=(a_{\alpha_1},\ldots,a_{\alpha_n})$, $\alpha_1<\ldots<\alpha_n$ such that $S\subseteq R({\bar a})^\mathrm{rcl}$. Then by choice of the $a_\alpha$ the chain $R\subseteq R(a_{\alpha_1})^\mathrm{rcl}\subseteq \ldots\subseteq R(a_{\alpha_1},\ldots,a_{\alpha_n})^\mathrm{rcl}$ witnesses that the value group of $R(a_{\alpha_1},\ldots,a_{\alpha_n})^\mathrm{rcl}$ has height $n$, where height stands for the number of convex subgroups of the value group. Since $\trdeg(S/R)=2$, there are $n-2$ elements $b_1,\ldots,b_{n-2}$ from $\{a_{\alpha_1},\ldots,a_{\alpha_n}\}$ that are algebraically independent over $S$. Since $S$ can be embedded into $R((t^\R))$ we know that $S$ has height $1$: Crucially we use here that any such embedding preserves the valuations because the natural valuation on $R((t^\R))$ again has the convex hull of $R$ in $R((t^\R))$ as its valuation ring. But now the chain $R\subseteq S\subseteq S(b_1)^\mathrm{rcl}\subseteq \ldots\subseteq S(b_1,\ldots,b_{n-2})^\mathrm{rcl}=R(a_{\alpha_1},\ldots,a_{\alpha_n})^\mathrm{rcl}$ witnesses that the value group of $R(a_{\alpha_1},\ldots,a_{\alpha_n})^\mathrm{rcl}$ has height at most $n-1$, which gives the desired contradiction. \end{proof} \begin{FACT}{Theorem}{CODFprimModels} Let $K$ be a differential and formally real field. If $K$ has a prime model extension $\hat K$ for $\CODF$, then $\hat K$ is algebraic over $K$. \end{FACT} \begin{proof} Suppose there is a prime model extension $\hat K$ of $K$ for $\CODF$ but $\hat K$ is not algebraic over $K$. Let $R$ be the algebraic closure of $K$ in $\hat K$. Then $R$ is a differential subfield of $\hat K$ and $\hat K$ is also a prime model of $R$ for $\CODF$: If $R\subseteq M\models \CODF$, then any $K$-embedding $\hat K\lra M$ must be the identity on $R$. Hence we may assume that $K$ is real closed all along. \par Choose real closed fields $M,N$ for $K$ as in \ref{RCFnoEmbeding}. By \ref{DLexistence} there are extensions of the derivation of $K$ to $M,N$ respectively such that $M,N$ furnished with these extensions are \CODF s. Since $\hat K$ can be embedded into $M$ and into $N$ by assumption, \ref{RCFnoEmbeding} implies that $\hat K$ must be of transcendence degree $\leq 1$ over $K$ and $K$ is Dedekind complete in $\hat K$. As $K\neq \hat K$, we know that $\trdeg(\hat K/K)=1$. \par Since $\hat K$ is a $\CODF$ it follows that $\hat K$ has a positive infinitesimal element $t$ with respect to $K$ such that $t'=1$ (in particular $t\notin K$). Then $\hat K$ is a differential subfield of $K((t^\Q))$ (endowed with the derivation extending the one on $K$ and satisfying $t'=1$). By \cite[end of 5.3]{LeSTre2024} we know that $t^{-1}$ has no integral in $K((t^\Q))$. This contradicts the fact that $t^{-1}$ has an integral in the \CODF\ $\hat K$. \end{proof} \begin{FACT}{Remark}{} The proofs of \ref{RCFnoEmbeding} and \ref{CODFprimModels} can be adapted to get the analogous statements about differential and formally p-adic fields and the class of p-adically closed differentially large fields. One possible task for future work is to extend \ref{CODFprimModels} (or rather \ref{RCFnoEmbeding}) to topological differential fields in the sense of \cite{GuzPoi2010}. We do not know if there is a version of \ref{CODFprimModels} outside of that context. For example, if $K$ is a subfield of a pseudo-finite field and $d$ is a derivation of $K$, it is unclear whether there is a prime model over $(K,d)$ in the class of differentially large and pseudo-finite fields (all of characteristic zero). \end{FACT} \par \renewcommand\href[3][]{}  \par \end{document}